\documentclass[3p]{elsarticle}

\usepackage[utf8]{inputenc}
\usepackage{graphicx,mathtools}
\usepackage{todonotes}
\usepackage{csvsimple}
 \usepackage[normalem]{ulem} 
\RequirePackage{amsmath,amsfonts,amssymb,stmaryrd}
\RequirePackage{mathtools}
\RequirePackage{bm}
\RequirePackage{doi}
\RequirePackage{multirow}
\RequirePackage{siunitx}
\RequirePackage{seqsplit}
\RequirePackage{xspace}
\RequirePackage{hhline}

\newcounter{para}[subsubsection]
\renewcommand\thepara{{\thesubsubsection}.\alph{para}}

\renewenvironment{paragraph}[1]%
{
\refstepcounter{para}
\vspace{1ex}
\textbf{\thepara\space-\space#1\space-\space}
}%

\graphicspath{{./figures/}}
\newcommand{\absperm}{\mathbb{K}}       
\newcommand{\normal}{\bm{n}}        
\newcommand{\porosity}{\phi}            
\newcommand{\head}{\ensuremath{h}}         
\newcommand{\source}{{q}}           
\renewcommand{\time}{{t}}           

\newcommand{\vecu}{\ensuremath{\bm{u}}}

\newcommand{\dirichlet}{{\footnotesize \head}} 
\newcommand{\neumann}{{\footnotesize u}}       

\newcommand{\K}{\absperm}

\newcommand{\globalDomain}{\Omega}
\newcommand{\originalDomain}{\Lambda}

\newcommand{\pLn}{{\partial\originalDomain_{\neumann}}}
\newcommand{\pLd}{{\partial\originalDomain_{\dirichlet}}}

\newcommand{\xmark}{\ensuremath{\times}}

\newcommand{\jump}[1]{\left\llbracket {#1} \right\rrbracket}

\newcommand{\aperture}{\varepsilon}

\newcommand{\UibTpfa}{\texttt{UiB-TPFA}\xspace}
\newcommand{\UibMpfa}{\texttt{UiB-MPFA}\xspace}
\newcommand{\UibRT}{\texttt{UiB-RT0}\xspace}
\newcommand{\UibMVEM}{\texttt{UiB-MVEM}\xspace}
\newcommand{\StuttTpfa}{\texttt{USTUTT-TPFA\_Circ}\xspace}
\newcommand{\StuttMpfa}{\texttt{USTUTT-MPFA}\xspace}
\newcommand{\Lanl}{\texttt{LANL-MFD}\xspace}

\newcommand{\Ncu}{\texttt{NCU\_TW-Hybrid\_FEM}\xspace}
\newcommand{\UniceVagC}{\texttt{UNICE\_UNIGE-VAG\_Cont}\xspace}
\newcommand{\UniceVagD}{\texttt{UNICE\_UNIGE-VAG\_Disc}\xspace}
\newcommand{\UniceHFVC}{\texttt{UNICE\_UNIGE-HFV\_Cont}\xspace}
\newcommand{\UniceHFVD}{\texttt{UNICE\_UNIGE-HFV\_Disc}\xspace}
\newcommand{\Ethz}{\texttt{ETHZ\_USI-FEM\_LM}\xspace}
\newcommand{\Unicamp}{\texttt{UNICAMP-Hybrid\_Hdiv}\xspace}
\newcommand{\Unil}{\texttt{UNIL\_USI-FE\_AMR\_AFC}\xspace}
\newcommand{\Inm}{\texttt{INM-EDFM}\xspace}
\newcommand{\Dtu}{\texttt{DTU-FEM\_COMSOL}\xspace}

\allowdisplaybreaks[1]
\title{Verification benchmarks for single-phase flow in three-dimensional fractured porous media}

\author[1]{Inga Berre} 
\author[2]{Wietse M. Boon}
\author[3]{Bernd Flemisch\corref{cor1}}
\author[1,4]{Alessio Fumagalli}
\author[3]{Dennis Gläser}
\author[1]{Eirik Keilegavlen}
\author[4]{Anna Scotti}
\author[1]{Ivar Stefansson}
\author[5,6]{Alexandru Tatomir}
\author[7]{Konstantin Brenner}
\author[8]{Samuel Burbulla}
\author[9]{Philippe Devloo}
\author[9]{Omar Duran}
\author[10]{Marco Favino}
\author[11]{Julian Hennicker}
\author[12,13]{I-Hsien Lee}
\author[14]{Konstantin Lipnikov}
\author[7]{Roland Masson}
\author[15]{Klaus Mosthaf}
\author[16]{Maria Giuseppina Chiara Nestola}
\author[12,13]{Chuen-Fa Ni}
\author[17]{Kirill Nikitin}
\author[18]{Philipp Schädle}
\author[14]{Daniil Svyatskiy}
\author[17]{Ruslan Yanbarisov}
\author[16]{Patrick Zulian}

\address[1]{Department of Mathematics, University of Bergen,
All\'{e}gaten 41, 5007 Bergen, Norway}
\address[2]{Department of Mathematics, KTH Royal Institute of Technology, Lindstedtsv\"agen 25, 11428 Stockholm, Sweden}
\address[3]{Department of Hydromechanics and Modelling of Hydrosystems,
University of Stuttgart, Pfaffenwaldring 61, 70569 Stuttgart, Germany}
\address[4]{Laboratory for Modeling and Scientific Computing MOX, Politecnico di Milano,
p.za Leonardo da Vinci 32, 20133 Milano, Italy}
\address[5]{Department of Applied Geology, Geosciences Center,
University of G\"ottingen, Goldschmidtstrasse 3, 37077 G\"ottingen, Germany}
\address[6]{Department of Earth Sciences, Uppsala University, Villav\"agen 16, S-75236 Uppsala, Sweden}
\address[7]{University of C{\^o}te d'Azur, CNRS, INRIA, LJAD, Nice, France}
\address[8]{Institute of Applied Analysis and Numerical Simulation, University of Stuttgart, Pfaffenwaldring 57, 70569 Stuttgart, Germany}
\address[9]{FEC-Universidade Estadual de Campinas, R. Josiah Willard Gibbs 85 - Cidade Universitária, Campinas-SP, Brazil, CEP 13083-839}
\address[10]{Institute of Earth Sciences, University of Lausanne, Building Geopolis, UNIL-Mouline, 1015 Lausanne, Switzerland}
\address[11]{Section de Mathématiques, Université de Genève, 2-4 rue du Lièvre, CP 64, 1211 Genève, Switzerland}
\address[12]{Graduate Institute of Applied Geology, National Central University, Taiwan}
\address[13]{Center for Environmental Studies, National Central University, Taiwan}
\address[14]{Los Alamos National Laboratory, New Mexico, USA}
\address[15]{Department of Environmental Engineering, Technical University of Denmark, Bygningstorvet,
Building 115, 2800 Kgs. Lyngby, Denmark}
\address[16]{Numerical Simulation in Science, Medicine and Engineering Group,
Institute of Computational Science, Università della Svizzera italiana.
Via G. Buffi 13, 6900 Lugano
Ticino, Switzerland}
\address[17]{Marchuk Institute of Numerical Mathematics of Russian Academy of Sciences, Moscow, Russia}
\address[18]{ETH Zürich, Geothermal Energy and Geofluids Group, Institute of
Geophysics, 8092 Zürich, Switzerland}
\cortext[cor1]{Corresponding author, \texttt{bernd@iws.uni-stuttgart.de}}

\begin{document}

\begin{abstract}
Flow in fractured porous media occurs in the earth's subsurface, in biological
tissues, and in man-made materials. Fractures have a dominating influence on flow processes, and the last decade has seen an extensive development of models and numerical methods that explicitly account for their presence. To support these developments, we present a portfolio of four benchmark cases for single-phase flow in three-dimensional fractured porous media. The cases are specifically designed to test the methods' capabilities in handling various complexities common to the geometrical structures of fracture networks. Based on an open call for participation, results obtained with 17 numerical methods were collected. This paper presents the underlying mathematical model, an overview of the features of the participating numerical methods, and their performance in solving the benchmark cases.  
\end{abstract}

\maketitle

\section{Introduction}
Flow in fractured porous media is characterized by an interaction between
the fractures and the surrounding porous medium, commonly referred to as the
matrix.
The strong influence of fracture network geometry on flow patterns has motivated the development of mathematical models and numerical methods that explicitly account for the geometry of fractures \cite{Berre:2018:FFP}. Considering flow both in the fractures and in the surrounding porous medium, these models are based on the conceptual discrete-fracture-matrix (DFM) representation of the fractured porous media.  

With the development of a wealth of simulation tools for flow in fractured porous media, a need for verification benchmarks for numerical methods has emerged. 
To accommodate this need, four research groups working in the field initiated a comparison study, which led to the presentation of a suite of two-dimensional benchmark tests and corresponding results for a range of numerical methods \cite{Flemisch:2018:BSF}. 
The methods were probed on test cases featuring known difficulties for numerical methods, including fracture intersections and combinations of blocking and conducting fractures.
The study exposed the relative strengths and weaknesses between the participating methods, both in terms of accuracy and computational cost.
After the publication of the results, these benchmark cases have been widely applied by the scientific community in testing numerical methods and new simulation tools \cite{arraras2018monolithic, budisa2019block,budisa2019mixed,fumagalli2019dual,koppel2018lagrange, koppel2019stabilized, odsaeter2019simple, schadle20193d, stefansson2018finite}. 

Based on the reception of the verification benchmarks \cite{Flemisch:2018:BSF} and the capabilities of three-dimensional modeling in the research community, the next phase in the work on verification benchmarks was launched with a call for participation \cite{berre2018call}. The purpose of this call was to extend the platform of verification benchmarks for numerical
methods to three-dimensional problems. In addition, the studies were extended to include simulations of linear tracer transport as a means to highlight additional nuances in the comparison of the calculated flow fields. The present paper discusses the results we received as answers to this call. 

The paper is organized as follows. In Section \ref{sec:process}, an overview of the participation process is given. In Section \ref{sec:model}, we describe the mathematical models for fluid flow and transport in fractured porous media. Section \ref{sec:discretization} briefly describes the participating numerical methods as well as the discretization of the transport problem. The four test cases are described in Section \ref{sec:cases}, with each description followed by a presentation and discussion of the corresponding results. Section \ref{sec:discussion} summarizes the discussion of the results, and Section \ref{sec:conclusion} provides concluding remarks.

\section{Benchmark Process}
\label{sec:process}

The publication of this verification benchmark study was laid out as a four-stage process: the development of benchmark cases, a call for participation, collection and synchronization of the results by the participants, and a final discussion and reporting. 

The process started with the participants of the benchmark study \cite{Flemisch:2018:BSF} developing four new test cases. These were designed to test the capabilities of numerical methods for DFM representations of flow in three-dimensional fracture networks. The design of each test case was led by the "benchmark case designers" listed in Section \ref{sec:cases}. An open call for participation was launched in September 2018 \cite{berre2018call}, followed by a dedicated mini-symposium at the SIAM Conference on Mathematical and Computational Issues in the Geosciences, March 2019, Houston. Researchers interested in participating in the benchmark followed a predefined registration procedure, were approved by the authors issuing the call, and were asked to sign a participation agreement. During this process, we received applications concerning 15 additional numerical methods, all of which were approved. Finally, the results of 12 of these methods were submitted and included in the study. 

The case descriptions presented in the call \cite{berre2018call} were accompanied by data in the form of geometry descriptions, existing simulation results, and plotting scripts,
all available in the Git repository \url{https://git.iws.uni-stuttgart.de/benchmarks/fracture-flow-3d.git}. 
This repository was reused in the fully transparent collection and synchronization phase.
During this phase, the results were uploaded and made available to all participants, and recomputations and adjustments were allowed until August 2019.
In the fourth phase, all participants contributed to the reporting of the results presented in Section \ref{sec:cases}. The last two phases were led by assigned "benchmark case coordinators". 
While access to the Git repository was restricted to the benchmark participants during the phase of collection and comparison of the results, all data have been made publicly available upon submission of this manuscript. In addition to the data and plotting scripts, five Jupyter notebooks are provided, four focusing on reproducing the figures encountered in Section \ref{sec:cases}, and one for facilitating the comparison of new results.

\section{Mathematical Models}
\label{sec:model}

We introduce two models for flow and transport in fractured media. First, the flow model is presented in the conventional equidimensional setting, allowing a natural introduction to the physical parameters. From this formulation, we derive the mixed-dimensional model through appropriate reduction of the equations. The mixed-dimensional model forms the focus of this study. Finally, we present the equi- and mixed-dimensional transport models.

\subsection{Equidimensional Flow Model}
\label{sec:equi_dimensional_flow_model}

We consider a steady-state, incompressible, single-phase flow through a porous medium
described by Darcy's law. With the imposition of mass conservation, the governing system of equations is given by
\begin{subequations}\label{eq:strong_dual}
\begin{gather}\label{eq:darcy}
	\begin{aligned}
 		\vecu + \K \nabla \head &= 0,  \\
 		\nabla \cdot \vecu &= \source,
    \end{aligned}
    \quad \text{in } \originalDomain.
\end{gather}
Here, $\vecu$ denotes the fluid velocity in \si{\metre\per\second}, $\K$ is hydraulic conductivity measured in \si{\metre\per\second}, $\head$ is hydraulic head measured in \si{\metre}, and $q$ represents a source/sink term measured in \si{\per\second}. The domain $\originalDomain \subset \mathbb{R}^3$ will be called the
equidimensional domain. The following boundary
conditions on the boundary $\partial \originalDomain$ of
$\originalDomain$ complete model \eqref{eq:darcy}:
\begin{gather} \label{eq:darcybc}
	\begin{aligned}
 		\head |_{\pLd}&= \overline{\head} & \quad \text{on } \pLd, \\
		\vecu\cdot\normal |_{\pLn} &= \overline{u} & \text{on } \pLn.
	\end{aligned}
\end{gather}
\end{subequations}
We assume $\partial \originalDomain = \pLd \cup \pLn$, $\pLd \cap \pLn =
\emptyset$, and $|\pLd| > 0$. 
In \eqref{eq:darcybc} $\cdot |_A$ is a suitable trace operator on $A \subset \partial \originalDomain$. $\overline{\head}$ is the hydraulic head imposed on the boundary $\pLd$, while
$\overline{u}$ is the prescribed Darcy velocity normal to the boundary $\pLn$ with respect to the outer unit normal vector $\normal$.

By substituting Darcy's law in the mass conservation, the dual problem \eqref{eq:strong_dual} can be recast in its primal formulation, given by
\begin{gather}\label{eq:strong_primal}
	\begin{aligned}
	    -\nabla \cdot \K \nabla \head &=\source &\quad \text{in } \originalDomain, \\
    	\head |_{\pLd} &=  \overline{\head} & \text{on } \pLd, \\
		- \K \nabla \head \cdot\normal|_{\pLn} &=  \overline{u} & \text{on } \pLn.
    \end{aligned}
\end{gather}
Problems \eqref{eq:strong_dual} and \eqref{eq:strong_primal} are equivalent. However, different
numerical schemes are based on either of the two formulations. Under regularity assumptions on $\originalDomain$ 
and the data, 
both problems admit a unique weak solution. We refer to \cite{Raviart1977,Brezzi1991,Roberts1991,Ern2004}
for more details. 

We assume that $\originalDomain$ contains several fractures, i.e., thin inclusions in the domain. The fracture walls are assumed to be planar with smooth boundaries. The fractures have two distinguishing features: (1) the thickness, which we measure by the aperture $\aperture$, is small compared to the extension of the fracture; and
the (2) hydraulic conductivity may differ significantly from that of the rest of $\originalDomain$. The latter implies that the fractures may have a significant impact on the flow in $\originalDomain$. 

We further make the assumption that the principal directions of the local hydraulic conductivity are aligned with the orientation of the fractures. In particular, the hydraulic conductivity in the matrix ($\K_3$), the fractures ($\K_2$), as well as in the intersections between two fractures ($\K_1$) and at the crossings of intersections ($\K_0$), can be decomposed in the following way:
\begin{align*}
	\K_3 &= K_3^{eq}, &
    \K_2 &= \begin{bmatrix}
    	K_2^{eq} & \begin{matrix} 0 \\ 0 \end{matrix} \\
        \begin{matrix} 0 & 0 \end{matrix} & \kappa_2^{eq}
	\end{bmatrix}, \\ 
    \K_1 &= \begin{bmatrix}
		K_1^{eq} & 0 & 0 \\
        0 & \kappa_1^{eq} & 0 \\
        0 & 0 & \kappa_1^{eq} \\
	\end{bmatrix}, & 
    \K_0 &= \begin{bmatrix}
		\kappa_0^{eq} & 0 & 0 \\
        0 & \kappa_0^{eq} & 0 \\
        0 & 0 & \kappa_0^{eq} \\
	\end{bmatrix}.
\end{align*} 
Here, $K_d^{eq}$ and $\kappa_d^{eq}$, for different values of $d$, denote the tangential and normal hydraulic conductivities, respectively. Thus, $K_d^{eq}$ is an elliptic $(d \times d)$-tensor field, whereas $\kappa_d^{eq}$ is a positive scalar field. The subscript $d$ indicates that the features will be represented by $d$-dimensional objects in the reduced model, as derived in the next section. The superscript $eq$, on the other hand, indicates that these quantities are related to the equidimensional model.

\subsection{Mixed-dimensional Flow Model}
\label{sec:mixed_dimensional_flow_model}

The small aperture of the fractures justifies a reduction of dimensionality 
to a representation where fractures and their intersections
are approximated by lower-dimensional objects. For more details on the derivation, we refer the reader to \cite{Alboin2000,Faille2002,Angot2003,Martin2005,DAngelo2011,Fumagalli2012g,Schwenck2015,Flemisch2016,Boon2018}.

Here, we use $\globalDomain$ to denote the mixed-dimensional decomposition of $\originalDomain$. First, let $\globalDomain$ contain a three-dimensional domain $\globalDomain_3$ that represents the (possibly unconnected) matrix. Furthermore, $\globalDomain$ contains up to three lower-dimensional, open
subdomains, namely, the union of fracture planes
$\globalDomain_{2}$, their intersection lines $\globalDomain_{1}$
and intersection points $\globalDomain_0$.
For compatibility, we assume that $\globalDomain_d \not\subset \globalDomain_{d'}$ for all $d' > d$. Finally, we introduce $\Gamma_{d} = \globalDomain_{d} \cap \partial \globalDomain_{d+1}$ as the set of $d$-interfaces
between neighboring subdomains of codimension one. Each interface is endowed with a normal unit vector $\normal$ pointing outward from $\globalDomain_{d+1}$. 

Remaining consistent with the notation convention above, data and unknowns will also be annotated with a subscript related to the dimension. As a first example, on a $d$-dimensional feature $\globalDomain_{d,i} \subseteq \globalDomain_d$ with counting index $i$, let $\aperture_{d,i}$ denote the cross-sectional volume, area, or length of the corresponding physical domain for $d = 0,...,2$, respectively. It has the unit of measure \si{\metre\tothe{3-d}} and is extended as nondimensional unity in $\globalDomain_3$. Moreover, we introduce for each $d$-feature with index $i$, a typical length $a_{d,i}$ such that $\aperture_{d,i} = a_{d,i}^{3 - d}$.
In the continuation, we will omit the subscript $i$ if no ambiguity arises.

We continue this subsection by presenting the reduced model associated with \eqref{eq:strong_dual} in the two-dimensional fractures $\globalDomain_2$ followed by its generalization for all $d = 0, ..., 3$.

\subsubsection{Two-dimensional Fracture Flow}

The variables in this formulation are the velocity $\vecu_3 = \vecu$ and hydraulic head $\head_3 = \head$ in the rock matrix $\globalDomain_3$, as well as the integrated tangential velocity $\vecu_2$ and average hydraulic head $\head_2$ in the fracture. These are given pointwise for $x \in \globalDomain_2$ by
\begin{gather*}
	\vecu_2(x) = \int_{\aperture_2(x)} \bm{u}_{\|}
    \quad \text{and} \quad
    \head_2(x) = \dfrac{1}{\aperture_2(x)} \int_{\aperture_2(x)} \head.
\end{gather*}
Here, $\bm{u}_{\|}$ denotes the components of $\bm{u}$ tangential to $\Omega_2$. The integrals are computed in the normal direction of the fracture, and thus, the corresponding units of measurement are \si{\metre\squared\per\second} and \si{\metre} for $\vecu_2$ and $h_2$, respectively.

We derive the reduced Darcy's law and the mass balance equation by averaging and integrating, respectively, over the direction normal to the fractures. Recall that the vector $\normal$ here refers to the normal unit vector oriented outward from $\globalDomain_3$.
\begin{subequations}\label{eq:strong_fracture_dual 2}
  \begin{gather}\label{eq:darcy_fracture_dual 2}
      \begin{aligned}
          \frac{1}{\aperture_2} \vecu_{2} + K_2^{eq} \nabla_2 \head_2 &= 0  \\
          \nabla_2 \cdot \vecu_2 - \jump{\vecu_3 \cdot \normal} &= \source_2
      \end{aligned}
      \quad \text{in } \globalDomain_2,
  \end{gather}
where $\nabla_2$ is the del-operator in the tangential directions and $\source_2$ is the integrated source term, i.e., $\source_2(s) = \int_{\aperture_2(s)} \source$. Note that we have assumed $K_2^{eq}$ to be constant in the direction normal to $\globalDomain_2$. The jump operator is defined as
$\jump{\vecu_3 \cdot \normal} |_{\globalDomain_d} = \sum (\vecu_3 \cdot \normal |_{\Gamma_2})$, thus representing the mass exchange between fracture and matrix. 
In particular, for each subdomain $\globalDomain_{2, i} \subseteq \globalDomain_2$, we sum over all flux contributions over sections of $\Gamma_{2}$ that coincide geometrically with $\globalDomain_{2, i}$.
These fluxes are assumed to satisfy the following Darcy-type law given by a finite difference between the hydraulic head in $\globalDomain_2$ and on $\partial\globalDomain_3$:
\begin{gather} \label{eq: Darcy normal 2}
	\vecu_3 \cdot \normal + \kappa_2^{eq} \frac{2}{a_d} (\head_2 - \head_3) = 0
    \quad \text{on } \Gamma_{2}.
\end{gather}
\end{subequations}
Note that to be mathematically precise, each term in this equation represents an appropriate trace or projection of the corresponding variable onto $\Gamma_2$.

\subsubsection{Generalized Flow Model}

Next, we generalize the equations described above to domains of all dimensions, thus including the intersection lines and points. For that purpose, we introduce the integrated velocity $\bm{u}_d$ for $d = 1$ and average hydraulic head $\head_d$ with $d = 0, 1$ given pointwise for $x \in \globalDomain_d$ by
\begin{gather*}
	\vecu_1(x) = \int_{\aperture_1(x)} \bm{u}_{\|}
    \quad \text{and} \quad
    \head_d(x) = \dfrac{1}{\aperture_d(x)} \int_{\aperture_d(x)} \head, \text{ for } d = 0, 1.
\end{gather*}
Again, $\bm{u}_{\|}$ denotes the components of $\bm{u}$ tangential to $\Omega_1$. 
The corresponding units of measurement are \si{\metre\cubed\per\second}  and \si{\metre} for $\vecu_1$ and $h_d$, respectively. The analogs of \eqref{eq:darcy_fracture_dual 2} on these lower-dimensional manifolds are then given by
  \begin{gather}\label{eq:darcy_fracture_dual d}
      \begin{aligned}
          \frac{1}{\aperture_1} \vecu_{1} + K_1^{eq} \nabla_1 \head_1 &= 0  \\
          \nabla_1 \cdot \vecu_1 - \jump{\vecu_2 \cdot \normal} &= \source_1
      \quad \text{in } \globalDomain_1, \\
      - \jump{\vecu_1 \cdot \normal} &= \source_0
      \quad \text{in } \globalDomain_0.
      \end{aligned}
  \end{gather}
Here, $\nabla_1$ denotes the del-operator, i.e., the derivative, in $\Omega_1$. Moreover, the linear jump operator $\jump{\cdot}$ is naturally generalized to
$\jump{\vecu_{d+1} \cdot \normal} |_{\globalDomain_d} = \sum (\vecu_{d+1} \cdot \normal |_{\Gamma_{d}})$, where we for each subdomain $\globalDomain_{d, i} \subseteq \globalDomain_d$ sum over all flux contributions over sections of $\Gamma_{d}$ that coincide geometrically with $\globalDomain_{d, i}$.  Finally, $q_1$ and $q_0$ correspond to the integrated source terms in the intersection lines and points, respectively.

Due to our choice of defining $\bm{u}_d$ as the integrated velocity, a scaling with $\aperture_{d + 1}$ appears in the equation governing the flux across $\Gamma_d$:
\begin{gather} \label{eq: Darcy normal d}
	\frac{1}{\epsilon_{d + 1}} \vecu_{d + 1} \cdot \normal + \kappa_d^{eq} \frac{2}{a_d} (\head_d - \head_{d + 1}) = 0
    \quad \text{on } \Gamma_{d}, \ d = 0, 1.
\end{gather}
Recalling that $\epsilon_3 = 1$, it now follows that the effective tangential and normal hydraulic conductivities are given by:
\begin{subequations}
	\begin{align}
      K_d &= \aperture_d K_d^{eq}, 
      &\text{in } \globalDomain_d, \ d &= 1, \dots, 3 \\
      \kappa_d &=  \aperture_{d + 1} \frac{2}{a_d} \kappa_d^{eq}, 
      &\text{on } \Gamma_d, \ d &= 0, \dots, 2.
	\end{align}
\end{subequations}
From these definitions, it is clear that the units of $K_d$ and $\kappa_d$ are \si{\metre\tothe{4-d}\per\second} and \si{\metre\tothe{2-d}\per\second}, respectively.

Collecting the above equations, we obtain the generalization of system \eqref{eq:strong_fracture_dual 2} to subdomains of all dimensions. The resulting system consists of Darcy's law in both tangential and normal directions followed by the mass conservation equations: 
\begin{subequations}\label{eq:strong_fracture_dual detailed}
	\begin{align}
 		\vecu_d + K_d \nabla_d \head_d &= 0, 
        &\text{in } &\globalDomain_d, \ d = 1, \dots, 3,\\
		\vecu_{d+1} \cdot \normal + \kappa_d (\head_d - \head_{d+1}) &= 0,
    	&\text{on } &\Gamma_d, \ d = 0, \dots, 2,\\
        \nabla_d \cdot \vecu_3 &= \source_3,
    	&\text{in } &\globalDomain_3, \\
 		\nabla_d \cdot \vecu_d - \jump{\vecu_{d+1} \cdot \normal} &= \source_d,
    	&\text{in } &\globalDomain_d, \ d = 1, 2,\\
 		- \jump{\vecu_1 \cdot \normal} &= \source_0,
    	&\text{in } &\globalDomain_0.
    \end{align}
\end{subequations}
The source term is given by $\source_3$ for the rock matrix and $\source_d(x) = \int_{\aperture_d(x)} \source$ measured in \si{\metre\tothe{3-d}\per\second}.

System \eqref{eq:strong_fracture_dual detailed} is then compactly described by:
\begin{subequations}\label{eq:strong_fracture_dual}
	\begin{align}
 		\vecu_d + K_d \nabla_d \head_d &= 0, 
        &\text{in } \globalDomain_d, \ d &= 1, \dots, 3, \label{eq:strong_Darcy_t} \\
		\vecu_{d+1} \cdot \normal + \kappa_d (\head_d - \head_{d+1}) &= 0,
    	&\text{on } \Gamma_d, \ d &= 0, \dots, 2, \label{eq:strong_Darcy_n} \\
 		\nabla_d \cdot \vecu_d - \jump{\vecu_{d+1} \cdot \normal} &= \source_d,
    	&\text{in } \globalDomain_d, \ d &= 0, \dots, 3, \label{eq:strong_massconv}
    \end{align}
\end{subequations}
in which the nonphysical $\vecu_4$ and $\vecu_0$ are understood as zero. 
The boundary conditions are inherited from the equidimensional model with the addition of a no-flux condition at embedded fracture endings:
\begin{subequations}\label{eq:strong_fracture_dual_BC}
	\begin{align}
        \head_d &= \overline{\head} 
        & \text{on } \partial \globalDomain_d \cap \pLd, \ d &= 0, \dots, 3,\\
		\vecu_d \cdot \normal &= \aperture_d \overline{u} 
        & \text{on } \partial \globalDomain_d \cap \pLn, \ d &= 1, \dots, 3, \\
        \vecu_d \cdot \normal &= 0 
        & \text{on } \partial \globalDomain_d \backslash (\Gamma_{d - 1} \cup \partial \originalDomain), \ d &= 1, \dots, 3.
    \end{align}
\end{subequations}

To finish the section, we present the primal formulation of the mixed-dimensional fracture flow model. Analogous to \eqref{eq:strong_primal}, this formulation is derived by substituting Darcy's laws \eqref{eq:strong_Darcy_t} and \eqref{eq:strong_Darcy_n} into the conservation equation \eqref{eq:strong_massconv}:
\begin{align}\label{eq:strong_fracture_primal}
	-\nabla_d \cdot K_d \nabla_d \head_d + \jump{\kappa_d (\head_d - \head_{d+1})} 
    &= \source_d,
   	&\text{in } \globalDomain_d, \ d &= 0, \dots, 3.
\end{align}
Again, we interpret the divergence term as zero if $d = 0$ and the jump term as zero if $d = 3$. The boundary conditions are given by
\begin{subequations}\label{eq:strong_fracture_primal_BC}
\begin{align}
        \head_d &= \overline{\head} 
        & \text{on } \partial \globalDomain_d \cap \pLd, \ d &= 0, \dots, 3,\\
		- K_d \nabla_d \head_d \cdot \normal &= \aperture_d \overline{u} 
        & \text{on } \partial \globalDomain_d \cap \pLn, \ d &= 1, \dots, 3, \\
        - K_d \nabla_d \head_d \cdot \normal &= 0 
        & \text{on } \partial \globalDomain_d \backslash (\Gamma_{d - 1} \cup \partial \originalDomain), \ d &= 1, \dots, 3.
\end{align}
\end{subequations}

Many discretization schemes presented in this study ignore flow in the one-dimensional fracture intersections and zero-dimensional intersections thereof. Although these correspond to discretizing a simpler model, this is perfectly in line with the proposed study.

\subsection{Equidimensional Transport Model}
\label{sec:equi_dimensional_transport_model}
We now consider a scalar quantity $c$ with the unit of measure \si{\metre\tothe{-3}}, which is transported through the porous medium subject to the velocity field resulting from the flow model presented in the previous sections. The purely advective transport of $c$ is described by the conservation equation:
\begin{equation}\label{eq:strong_transport}
	\porosity \frac{\partial c}{\partial \time} + \nabla \cdot \left( c \vecu \right) = \source_c	\quad \text{in } \originalDomain,
\end{equation}
where $\porosity$ is the porosity of the medium and $\source_c$ is a source/sink term for $c$ given in \si{\metre\tothe{-3}\per\second}. We define Dirichlet boundary conditions on those boundary segments where inflow occurs, i.e.,
\begin{equation} \label{eq:transportbc}
 		c |_{\partial\originalDomain_{\footnotesize c}}
            = \overline{c} \quad \text{on } \partial\originalDomain_{\footnotesize c}, \quad
              \partial\originalDomain_{\footnotesize c} =
                 \{ x \in \partial\originalDomain: \,\, \vecu \cdot \normal < 0 \},
\end{equation}
with $\overline{c}$ being the value for $c$ prescribed on the boundary $\partial\originalDomain_{\footnotesize c}$.

\subsection{Mixed-dimensional Transport Model}
\label{sec:mixed_dimensional_transport_model}
Analogous to Section \ref{sec:mixed_dimensional_flow_model}, we choose the average value for $c$ as the primary variable, which is defined as $c_3 = c$ in $\globalDomain_3$ and for the lower dimensional objects (with $d \le 2$) as
\begin{equation*}
c_d(s) = \dfrac{1}{\aperture_d(s)} \int_{\aperture_d(s)} c.
\end{equation*}
Following the derivation of the mixed-dimensional flow model presented in Section \ref{sec:mixed_dimensional_flow_model}, the resulting mixed-dimensional transport model reads as:
\begin{equation}\label{eq:strong_transport_dual}
\aperture_d \porosity_d \frac{\partial c_d}{\partial \time}
		 + \nabla_d \cdot \left( c_d \vecu_d \right)
		 - \jump{\tilde{c}_{d+1} \left( \vecu_{d+1} \cdot \normal \right)}
		 = \source_{c,d}	\quad \text{in } \globalDomain_d, 
           \quad d = 0, \dots, 3.
\end{equation}
Note that for $d=0$, the divergence term is void. 
Here, the porosity is simply $\porosity_d = \porosity^{eq}$, with units of measure \si{\metre\tothe{-3}}, and $\tilde{c}_{d+1}$ is evaluated on the basis of a first-order upwind scheme, i.e.,
\begin{equation} \label{eq:upwind}
  \tilde{c}_{d+1} =
        \begin{cases}
            c_{d+1} & \text{if} \quad \vecu_{d+1} \cdot \normal |_{\Gamma_{d}} > 0\\
            c_d   & \text{if} \quad \vecu_{d+1} \cdot \normal |_{\Gamma_{d}} < 0.
        \end{cases}
\end{equation}
As in the flow model, the jump operator represents the sum of the fluxes over all contributions defined on sections of $\Gamma_{d}$ that coincide geometrically with $\globalDomain_{d, i}$.

\section{Discretization Methods} \label{sec:discretization}
The intent of this benchmark study is to quantitatively evaluate different
discretization schemes for the mixed-dimensional flow models \eqref{eq:strong_fracture_dual}-\eqref{eq:strong_fracture_primal_BC}.
As a means of evaluation, the velocities were inserted into a standard cell-centered, first-order upwind scheme for the transport equations \eqref{eq:strong_transport_dual}. The temporal discretization is given by the implicit Euler method with a fixed time-step prescribed for each test case.
The main properties of the discretization methods covered by the benchmark are summarized in Tables \ref{tab:scheme_properties_general} and \ref{tab:scheme_properties_numerical}, which also contain references for further details. 
The majority of the methods followed the mixed-dimensional flow model and the specified transport discretization, with the following exceptions:

The schemes \Ncu and \Dtu describe the flow along the fractures by additional terms defined on the fracture surfaces.
This effectively adds connectivity between the degrees of freedom located on fractures without introducing additional degrees of freedom.
This means that these schemes do not solve the mass balances \eqref{eq:strong_massconv} for $d < 3$.
Moreover, this approach implies continuity of the hydraulic head across the fractures and therefore replaces the coupling condition \eqref{eq:strong_Darcy_n}.
Other schemes participating in this study also assume continuity of the hydraulic head across the fractures, and a complete overview is given in Table \ref{tab:scheme_properties_numerical}.

The scheme \Unil is an equidimensional approach, meaning that the fractures, their intersections, and intersections of intersections are discretized with three-dimensional elements using locally refined grids.
Therefore, the lower-dimensional mass balances \eqref{eq:strong_massconv} for $d < 3$ and the coupling conditions \eqref{eq:strong_Darcy_n} are not relevant for this scheme.

Finally, the schemes \Ethz and \Unil do not use a first-order upwind scheme but apply an algebraic flux correction technique for the stabilization of a finite element discretization of the transport model \cite{kuzmin2012flux}. Such stabilization techniques provide a similar discretization as the given upwind scheme.

\begin{table}[hbt]
    \centering
    \begin{tabular}{|c c c c c|}
        \hline
        Acronym & References & Open source code & Run scripts & Test cases \\ \hhline{|=====|}        \multicolumn{5}{|c|}{Two-point flux approximation} \\ 
        \UibTpfa &  \cite{keilegavlen2019porepy,Nordbotten2018} & \checkmark \cite{porepySourceCode} & \cite{porepyRunScripts} & 1-4 \\
        \hline
        \multicolumn{5}{|c|}{Multi-point flux approximation} \\ 
        \UibMpfa &  \cite{keilegavlen2019porepy, Nordbotten2018} & \checkmark \cite{porepySourceCode} & \cite{porepyRunScripts} & 1-4 \\
        \hline
        \multicolumn{5}{|c|}{Lowest order mixed virtual element method} \\ 
        \UibMVEM & \cite{keilegavlen2019porepy, Nordbotten2018} & \checkmark \cite{porepySourceCode} & \cite{porepyRunScripts} & 1-4 \\ \hline
        \multicolumn{5}{|c|}{Lowest order Raviart-Thomas mixed finite elements} \\ 
        \UibRT & \cite{keilegavlen2019porepy, Nordbotten2018, Boon2018} & \checkmark \cite{porepySourceCode} & \cite{porepyRunScripts} & 1-4 \\ \hline
        \multicolumn{5}{|c|}{Multi-point flux approximation} \\ 
        \StuttMpfa &  \cite{flemisch2011dumux} & \checkmark \cite{dumuxSourceCode} & \cite{dumuxRunScripts} & 1-4 \\ \hline
        \multicolumn{5}{|c|}{Two-point flux approximation} \\ 
        \StuttTpfa &  \cite{flemisch2011dumux} & \checkmark \cite{dumuxSourceCode} & \cite{dumuxRunScripts} & 1-4 \\ \hline
        \multicolumn{5}{|c|}{Mimetic Finite Differences} \\ 
        \Lanl &  \cite{lipnikov2014lanlmfd} & \checkmark \cite{amanziSourceCode} &  & 1-4 \\ \hline
        
        \multicolumn{5}{|c|}{Hybrid finite element method} \\ 
        \Ncu &  \cite{lee2015ncutwHybridFEM, Lee:2019:SMF} &  &  & 1 \\ \hline
        \multicolumn{5}{|c|}{Vertex Approximate Gradient continuous hydraulic head} \\ 
        \UniceVagC & \cite{brenner2016cont} &  &  & 1-4 \\ \hline
        \multicolumn{5}{|c|}{Hybrid Finite Volumes continuous hydraulic head} \\ 
        \UniceHFVC & \cite{brenner2016cont} &  &  & 1-4 \\ \hline
        \multicolumn{5}{|c|}{Vertex Approximate Gradient discontinuous hydraulic head} \\ 
        \UniceVagD & \cite{brenner2016disc} &  &  & 1-4 \\ \hline
        \multicolumn{5}{|c|}{Hybrid Finite Volumes discontinuous hydraulic head} \\ 
        \UniceHFVD & \cite{brenner2016disc} &  &  & 1-4 \\ \hline
        \multicolumn{5}{|c|}{Lagrange multiplier - L2-projection finite elements} \\ 
        \Ethz & \cite{schadle20193d,koeppel2019femlm,krause2016femlm} & \checkmark & \cite{utopiagit} & 1-4 \\ \hline
        \multicolumn{5}{|c|}{Hybrid H(div)} \\ 
        \Unicamp & \cite{Devloo2019,Duran2019}& \checkmark & \cite{unicampRunScripts} & 1-4 \\ \hline
        
         \multicolumn{5}{|c|}{Flux-corrected finite element method and adaptive mesh refinement} \\ 
        \Unil & \cite{favinofully,kuzmin2012flux} &  &  &  1-3 \\ \hline
        \multicolumn{5}{|c|}{Embedded discrete fracture method} \\ 
        \Inm & \cite{ny-jcam2020} & \xmark & \xmark & 1,3 \\ \hline
        \multicolumn{5}{|c|}{First-order Lagrangian finite elements (COMSOL)} \\ 
        \Dtu & \cite{comsol} & \xmark & \cite{thisGitRepo} & 1-4 \\ \hline
        
    \end{tabular}
    \caption{Names, acronyms, references and test cases covered for all participating discretization methods.} 
    \label{tab:scheme_properties_general}
\end{table}

\begin{table}[hbt]
    \centering
    \begin{tabular}{|c p{26mm} p{18mm} p{20mm} p{24mm} p{18mm}|}
        \hline
         {Acronym} &  Degrees of freedom &  Local mass conservation & Allows $\head$ discontinuity & Conformity & subdomain dimensions \\ \hhline{|======|}
         \UibTpfa & \head (elem), $\lambda$ (mortar flux) & \checkmark & \checkmark & geometrically & 0-3\\ \hline
         \UibMpfa & \head (elem), $\lambda$ (mortar flux) & \checkmark & \checkmark & geometrically & 0-3\\ \hline
         \UibMVEM & \head (elem), $\vecu$ (faces), $\lambda$ (mortar flux) & \checkmark & \checkmark & geometrically & 0-3 \\ \hline
         \UibRT & \head (elem), $\vecu$ (faces), $\lambda$ (mortar flux) & \checkmark & \checkmark & geometrically & 0-3 \\ \hline
         \StuttMpfa &  \head (elem) & \checkmark & \checkmark & fully & 2-3\\ \hline
         \StuttTpfa & \head (elem) & \checkmark & \checkmark & fully & 2-3\\ \hline
        \Lanl & \head (faces) & \checkmark & \checkmark & fully & 2-3\\ \hline
        \Ncu & \head, \vecu (nodes) & \checkmark & \xmark & fully & 2-3 \\ \hline
        \UniceVagC &  \head (nodes), (fracture faces) &  \checkmark &  \xmark & conforming & 2-3 \\\hline
        \UniceVagD & \head (nodes), (fracture faces) & \checkmark & \checkmark & conforming & 2-3 \\ \hline
        \UniceHFVC & \head (faces), (fracture edges)  & \checkmark & \xmark & conforming & 2-3\\ \hline
        \UniceHFVD & \head(faces), (fracture edges) &  \checkmark &  \checkmark  & conforming & 2-3\\ \hline
        \Ethz & \head (nodes) $\lambda$ (nodes) & \xmark & \xmark & none & 2-3 \\\hline
        \Unicamp  & \head, \vecu (elem), $\lambda$ (faces) & \checkmark & \checkmark & geometrically & 0-3 \\\hline
        \Unil & \head (nodes) &  \checkmark &  \xmark & 
        not applicable
        & equi-dim.\\\hline
        \Inm & \head (elem) & \checkmark & \xmark & none & 2-3 \\\hline
        \Dtu & \head (nodes) & \checkmark & \xmark & fully & 2-3 \\ 
        \hline
    \end{tabular}
    \caption{Numerical properties for the discretization methods. An entry in the column ``conforming'' can be  ``fully'' if each fracture element needs to coincide with a facet shared by two neighboring matrix elements, ``geometrically'' if each fracture needs to be a union of element facets from each of the two neighboring matrix subdomain meshes, or ``none'' if fracture and matrix meshes can be completely independent of each other.}
    \label{tab:scheme_properties_numerical}
\end{table}


\section{Benchmark Cases and Results}
\label{sec:cases}
In this section, we present the benchmark cases and compare the submitted results. For each case, the 
hydraulic head and tracer concentration are compared using several predefined
macroscopic metrics. In Subsection \ref{subsec:single_fracture}, a benchmark case containing a single fracture problem is
considered. Subsection \ref{subsec:regular} presents a benchmark based on a synthetic 
network composed of nine, regularly arranged fractures. The benchmark case in Subsection \ref{subsec:small_features} considers
the geometrically challenging case of almost intersecting fractures, fractures with small intersections, and
other features that a fracture network may exhibit. Finally, in Subsection
\ref{subsec:field_network}, we study a case with 52 fractures selected from a real network.

\subsection{Case 1: Single Fracture}\label{subsec:single_fracture}

\noindent
\textbf{Benchmark case designers:} D. Gläser and A. Tatomir\\
\textbf{Benchmark case coordinators:} B. Flemisch and A. Tatomir 

\subsubsection{Description}

\begin{figure}[hbtp]
	\centering
    \includegraphics[width=0.8\textwidth]{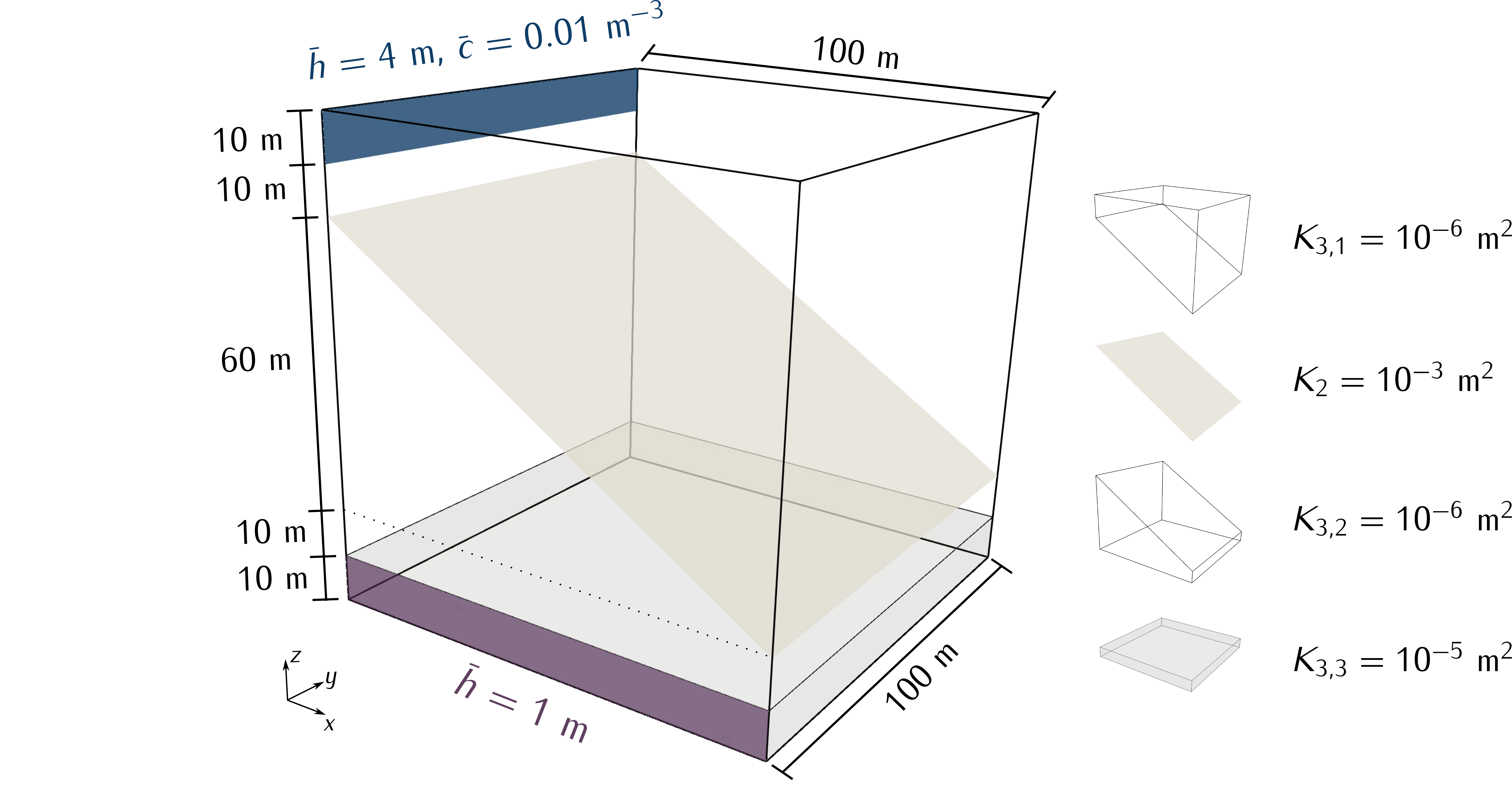}
    \caption{Conceptual model and geometrical description of the domain for Case 1 of Subsection \ref{subsec:single_fracture}.}
    \label{fig:problem1}
\end{figure}

Figure \ref{fig:problem1} illustrates the first benchmark case, with 
a geometry that is slightly modified from works \cite{Zielke:1991:DMT} and \cite{Barlag:1998:AMM}. 
The domain $\Omega$ is a cube-shaped region (\SI{0}{\metre}, \SI{100}{\metre}) $\times$ (\SI{0}{\metre}, \SI{100}{\metre}) $\times$ (\SI{0}{\metre}, \SI{100}{\metre}) which is crossed by a planar fracture, $\globalDomain_2$, with a thickness of \SI{1e-2}{\metre}. The matrix domain consists of subdomains $\globalDomain_{3,1}$ above the fracture and $\globalDomain_{3,2}$ and $\globalDomain_{3,3}$ below. The subdomain $\globalDomain_{3,3}$ represents a heterogeneity within the rock matrix. Inflow into the system occurs through a narrow band defined by $\{\SI{0}{\metre}\}\times (\SI{0}{\metre}, \SI{100}{\metre}) \times (\SI{90}{\metre}, \SI{100}{\metre}) $. Similarly, the outlet is a narrow band defined by $(\SI{0}{\metre}, \SI{100}{\metre}) \times \{\SI{0}{\metre}\}\times (\SI{0}{\metre}, \SI{10}{\metre})$. 

At the inlet and outlet bands, we impose the hydraulic head $h_{in} = \SI{4}{\metre}$ and $h_{out} = \SI{1}{\metre}$  respectively, and $c_{in} = \SI{1e-2}{\metre\tothe{-3}}$ is set at the inlet for the transport problem. All remaining parts of the boundary are assigned no-flow conditions. The parameters for conductivity, porosity, and aperture are listed in Table \ref{tab:parameters_case1} together with the overall simulation time and time-step size.

\begin{table}[hbt]
\begin{center}
\begin{tabular}{|l|ll|}\hline
Matrix hydraulic conductivity $K_{3,1}$, $K_{3,2}$ & \num{1e-6}$\bm{I}$ & \si{\metre\per\second} \\
Matrix hydraulic conductivity $K_{3,3}$         & \num{1e-5}$\bm{I}$ & \si{\metre\per\second} \\
Fracture effective tangential hydraulic conductivity $K_2$ & \num{1e-3}$\bm{I}$ & \si{\metre^2\per\second} \\
Fracture effective normal hydraulic conductivity $\kappa_2$ & \num{20} &\si{\per\second} \\
Matrix porosity $\phi_{3,1}, \phi_{3,2}$         & \num{2e-1}                              &\\
Matrix porosity $\phi_{3,3}$                   & \num{2.5e-1}                             & \\
Fracture porosity $\phi_2$                   & \num{4e-1}                              & \\
Fracture cross-sectional length $\epsilon_2$                        & \num{1e-2}                      & \si{\metre} \\
Total simulation time & \num{1e9} & \si{\second} \\
Time-step $\Delta t$  & \num{1e7} & \si{\second} \\ \hline
\end{tabular}
\end{center}
\caption{Parameters used in Case 1 of Subsection \ref{subsec:single_fracture}.}
\label{tab:parameters_case1}
\end{table}


\subsubsection{Results} 

Three different simulations were carried out with approximately \num{1}k, \num{10}k and \num{100}k cells for the 3d domain. The precise number of cells and degrees of freedom for each method are listed in Table \ref{tab:case1_comp_cost} and will be discussed in Subsection \ref{case1:g}. We compare the methods on the basis of computed pressure head and concentration, plotted along prescribed lines. The first comparison, represented in Subsection \ref{case1:a}, depicts the hydraulic head along a line crossing the 3d matrix domain, while the solutions reported in \ref{case1:b} and \ref{case1:c} visualize the matrix and fracture concentration along lines at the final simulation time. Plots in Subsection \ref{case1:d} and \ref{case1:e} depict integrated matrix and fracture concentrations over time, respectively. Finally, we compare concentration fluxes across the outlet over time in Subsection \ref{case1:f}.

\begin{figure}[hbtp]
	\centering
	\includegraphics[width=0.95\textwidth]{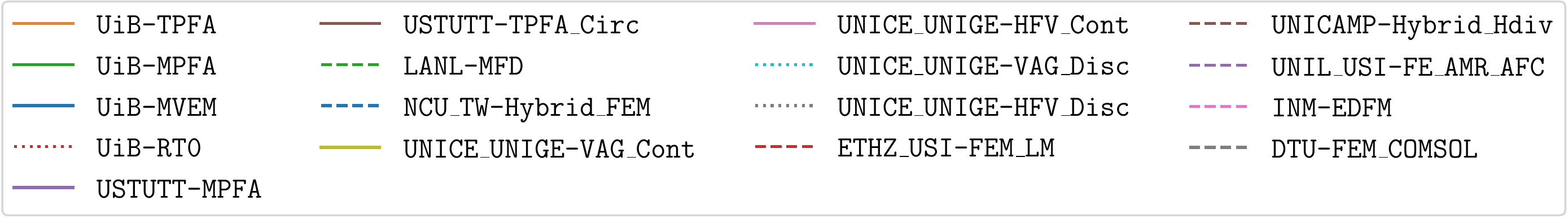}\\
	\includegraphics[width=0.95\textwidth]{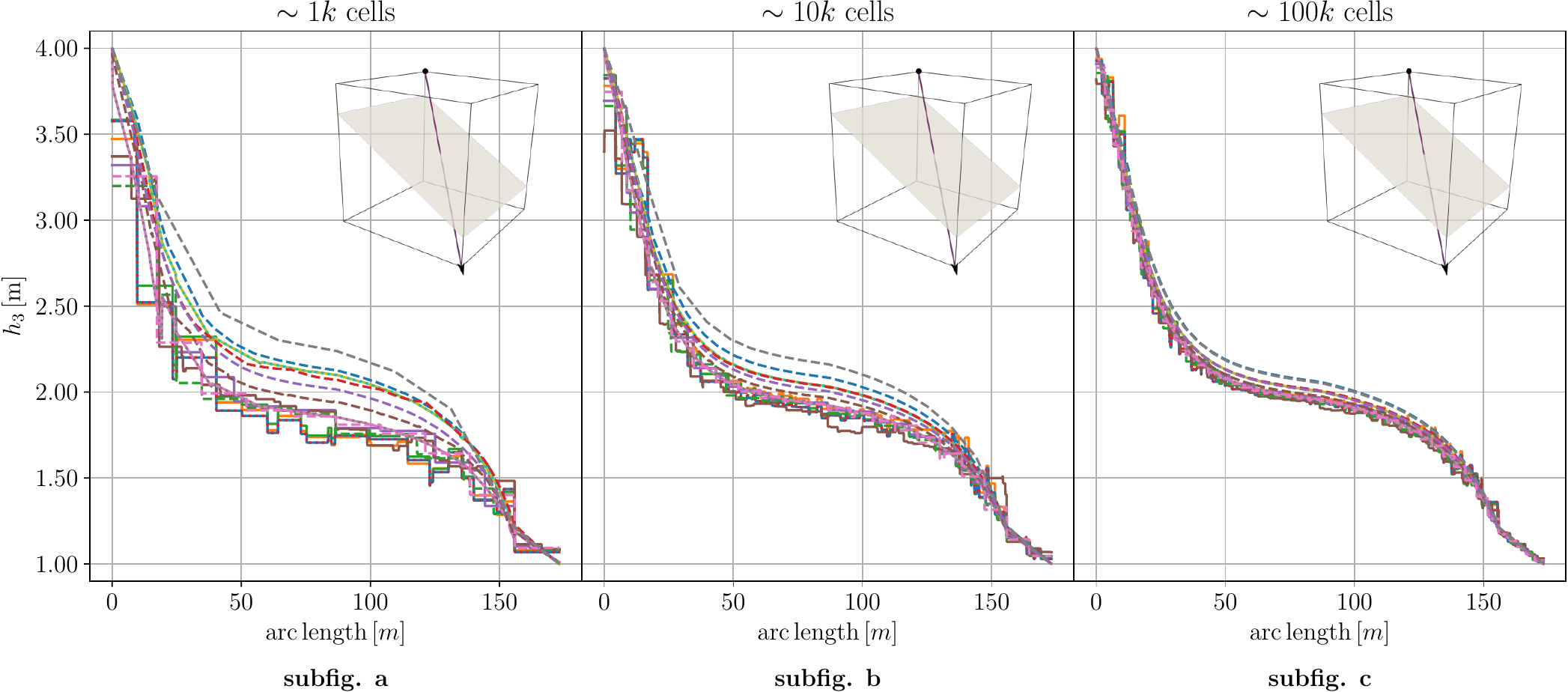}\\
	\includegraphics[width=0.95\textwidth]{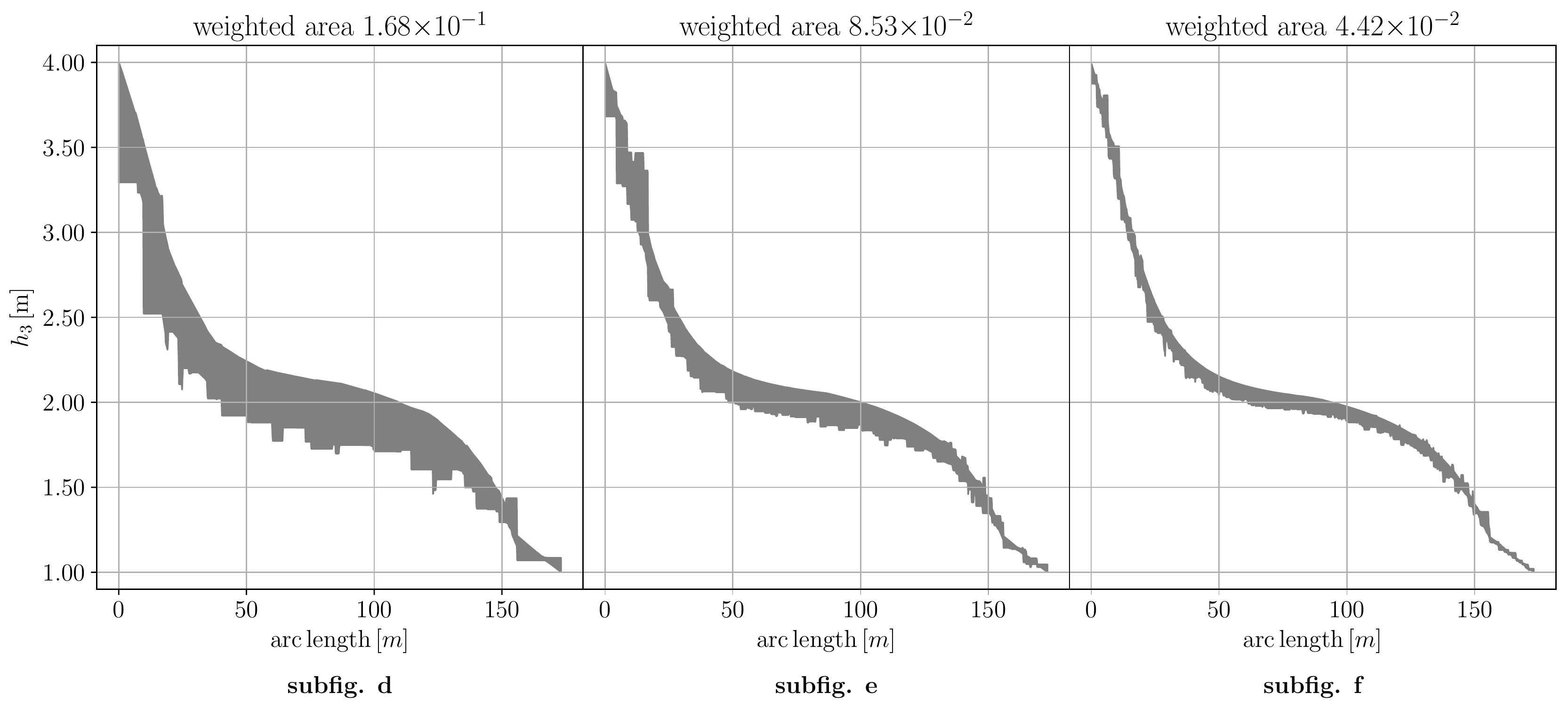}
    \caption{Case 1 of Subsection \ref{subsec:single_fracture}. On the top, the hydraulic head $\head_3$ in the matrix over the line $\left( \SI{0}{\metre}, \SI{100}{\metre}, \SI{100}{\metre} \right)$ - $\left( \SI{100}{\metre}, \SI{0}{\metre}, \SI{0}{\metre} \right)$ for three refinements (coarse to fine). On the bottom, the area between the 10th and 90th percentiles for three refinements (coarse to fine) and data. Results of Subsection \ref{case1:a}.}
    \label{fig:case1_pol_p_matrix}
\end{figure}

\paragraph{Hydraulic Head Over Line}\label{case1:a}
Figure \ref{fig:case1_pol_p_matrix} depicts the hydraulic head $\head_3$ in the matrix along the line $\left( \SI{0}{\metre}, \SI{100}{\metre}, \SI{100}{\metre} \right)$-$\left( \SI{100}{\metre}, \SI{0}{\metre}, \SI{0}{\metre} \right)$. Each plot corresponds to one of the three refinement levels.

At the coarsest level of around 1000 cells, all methods already show reasonable agreement. As expected, differences between the methods decrease with increasing refinement level. We remark that two classes of methods can be distinguished in these plots. First, the methods that use cellwise constant values exhibit staircase-like patterns. On the other hand, methods using nodal values are interpolated within each cell and yield a smoother appearance.

To quantify the differences between the participating methods and their convergence behavior over all refinement levels, we calculate and visualize the spread of the associated data sets. 
For that purpose, the solution values are evaluated at 1000 evenly distributed points along the considered line. At each such point, the mean as well as the 10th and 90th percentiles are determined. Each plot in the bottom row of Figure \ref{fig:case1_pol_p_matrix} visualizes the area between the 10th and 90th percentiles over the evaluation points. The number in the picture title corresponds to that area divided by the area under the mean curve. Convergence between the methods can clearly be observed.

\begin{figure}[hbtp]
	\centering
	\includegraphics[width=0.95\textwidth]{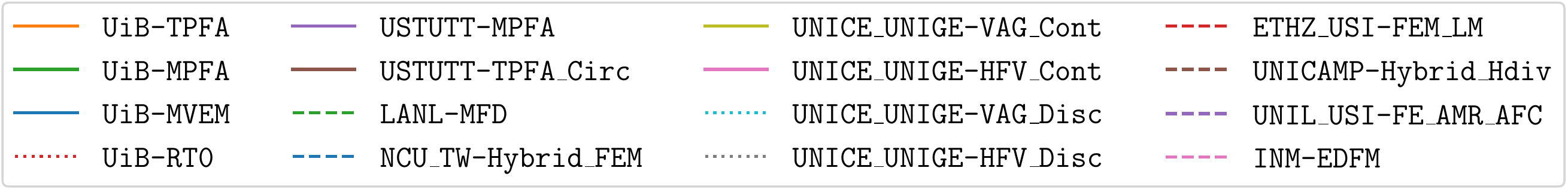}\\
	\includegraphics[width=0.95\textwidth]{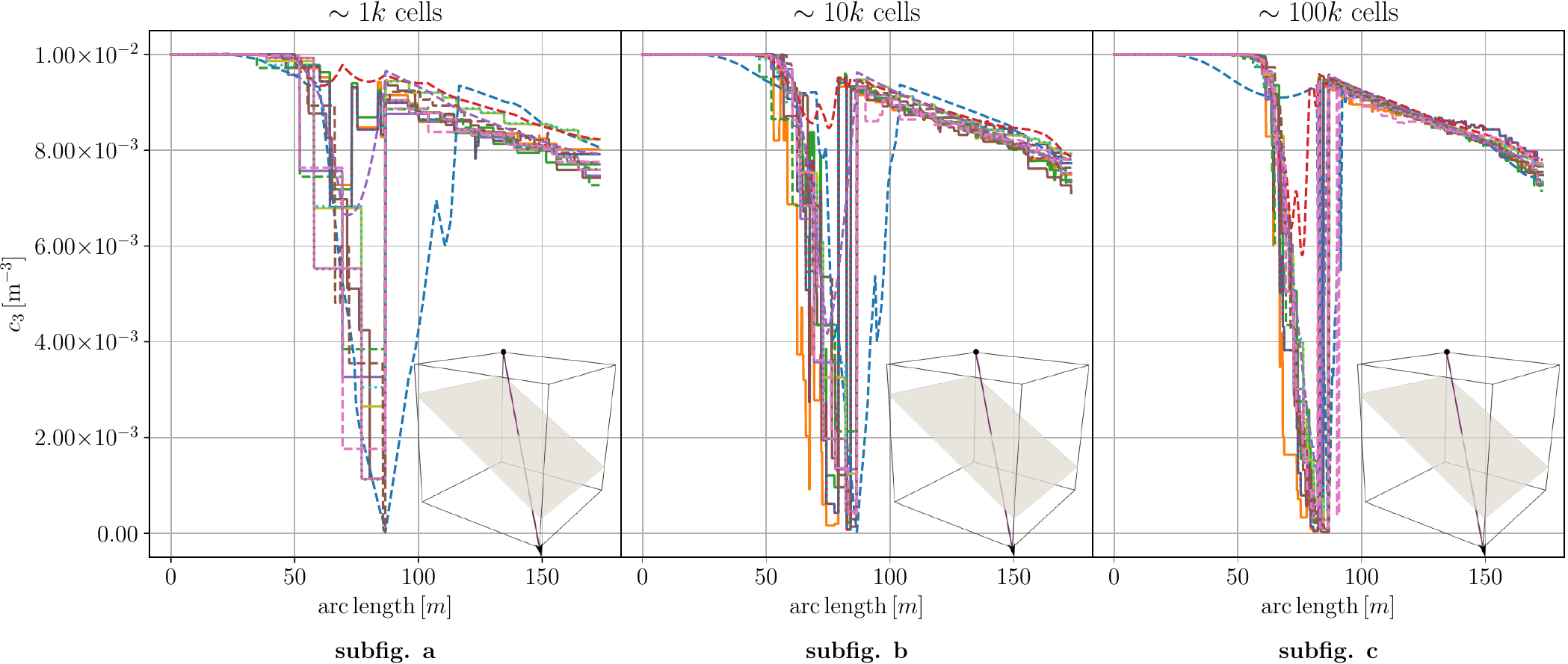}\\
	\includegraphics[width=0.95\textwidth]{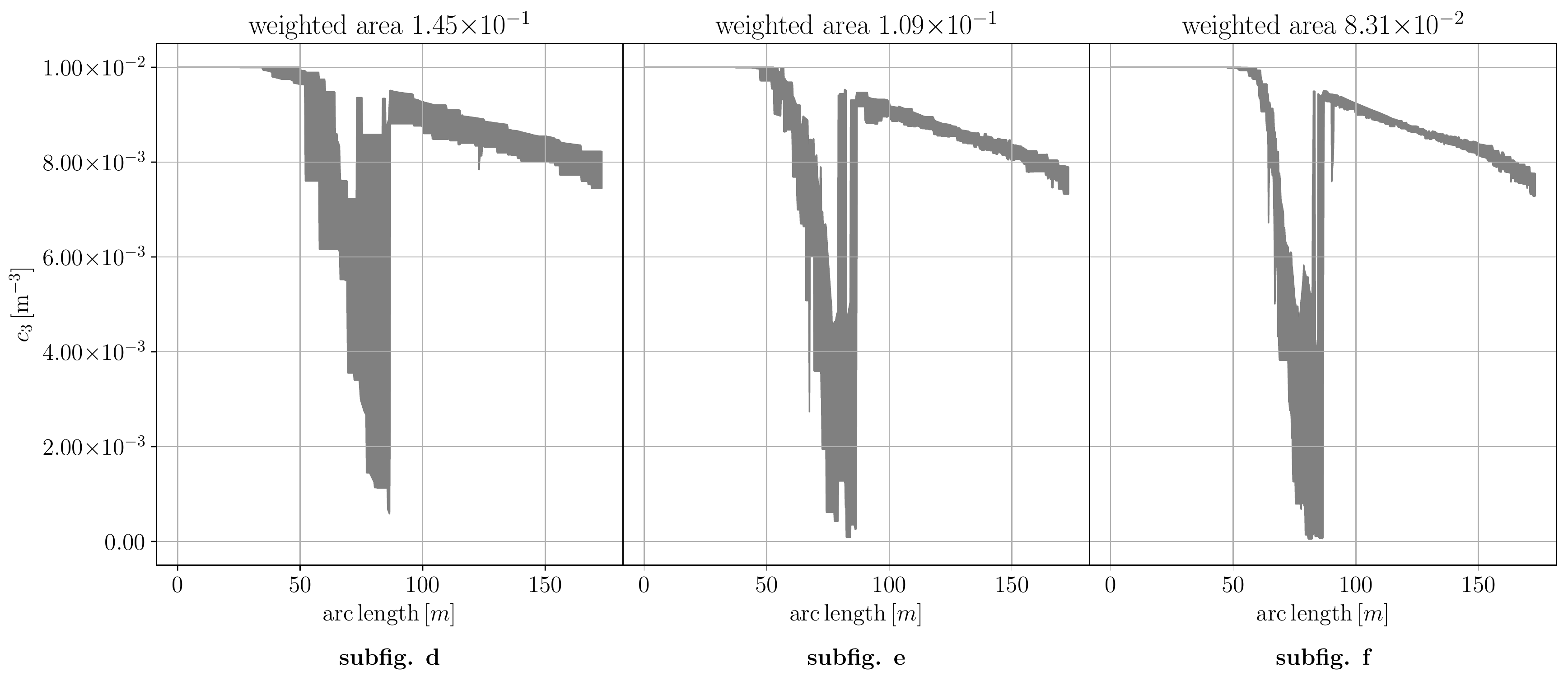}
    \caption{Case 1 of Subsection \ref{subsec:single_fracture}. On the top, concentration $c_3$ in the matrix, at the final simulation time, along the line $\left( \SI{0}{\metre}, \SI{100}{\metre}, \SI{100}{\metre} \right)$ - $\left( \SI{100}{\metre}, \SI{0}{\metre}, \SI{0}{\metre} \right)$ for three refinements (coarse to fine). On the bottom, area between the 10th and 90th percentiles for three refinements (coarse to fine) and data. Results of Subsection \ref{case1:b}.}
    \label{fig:case1_pol_c_matrix}
\end{figure}

\paragraph{Matrix Concentration Over Line}\label{case1:b}
The pictures at the top of Figure \ref{fig:case1_pol_c_matrix} illustrate the concentration $c_3$ in the matrix at the final simulation time along the line $\left( \SI{0}{\metre}, \SI{100}{\metre}, \SI{100}{\metre} \right)$-$\left( \SI{100}{\metre}, \SI{0}{\metre}, \SI{0}{\metre} \right)$, again for the different refinement levels.
We observe a similar behavior to that in \ref{case1:a} in the sense that the differences between most of the methods decrease with increasing refinement level. However, two methods show more pronounced deviations from the rest: \Ethz exhibits oscillations that can be attributed to the fact that the employed algebraic flux correction stabilization scheme does not suppress all spurious oscillations. The \Ncu does not capture the curve behavior at all. The obviously larger spread in the results is visualized more explicitly in the bottom row of Figure \ref{fig:case1_pol_c_matrix}. As a result, the convergence is much slower compared to Subsection \ref{case1:a}.

\begin{figure}[hbtp]
	\centering
	\includegraphics[width=0.95\textwidth]{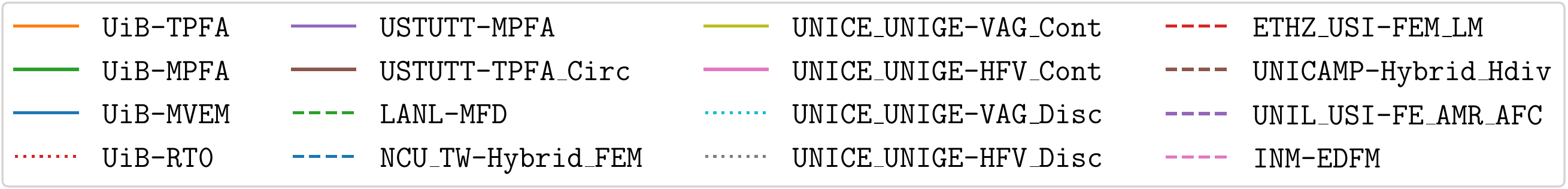}\\
	\includegraphics[width=0.95\textwidth]{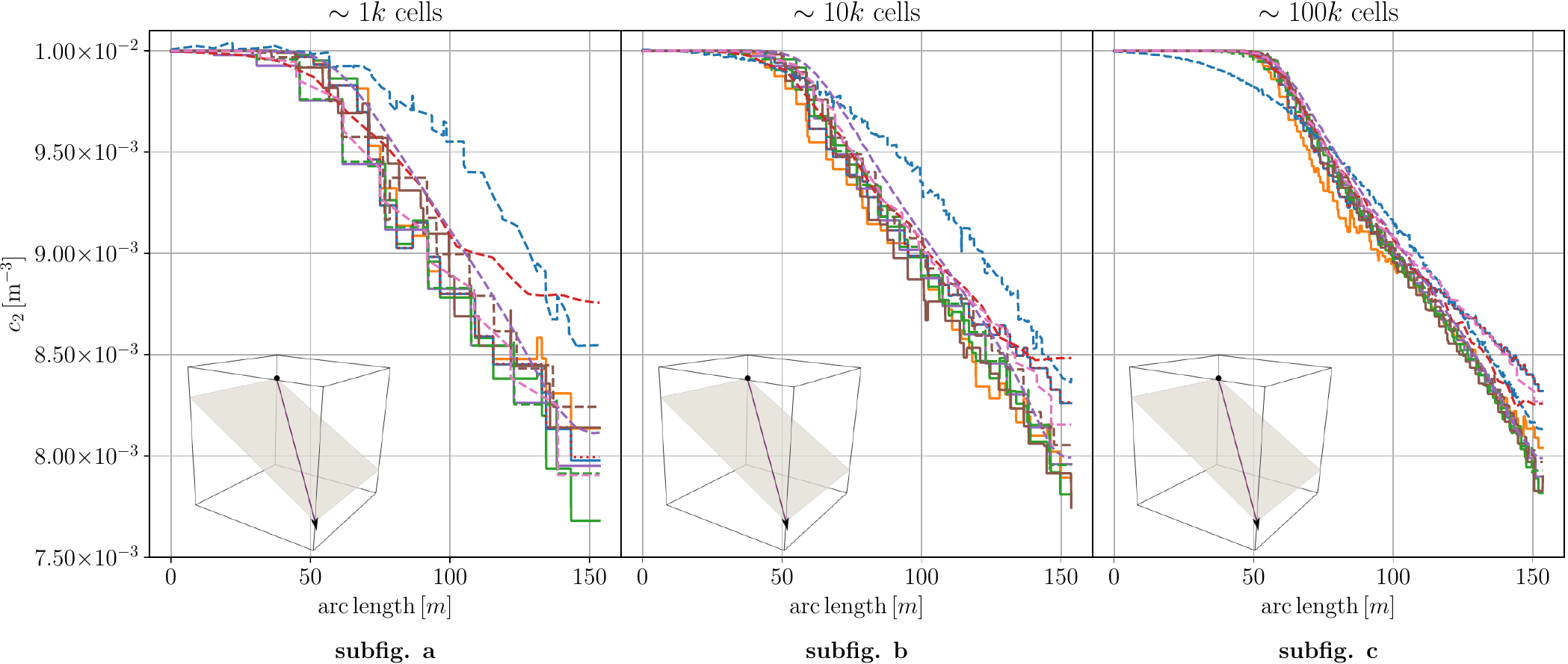}\\
	\includegraphics[width=0.95\textwidth]{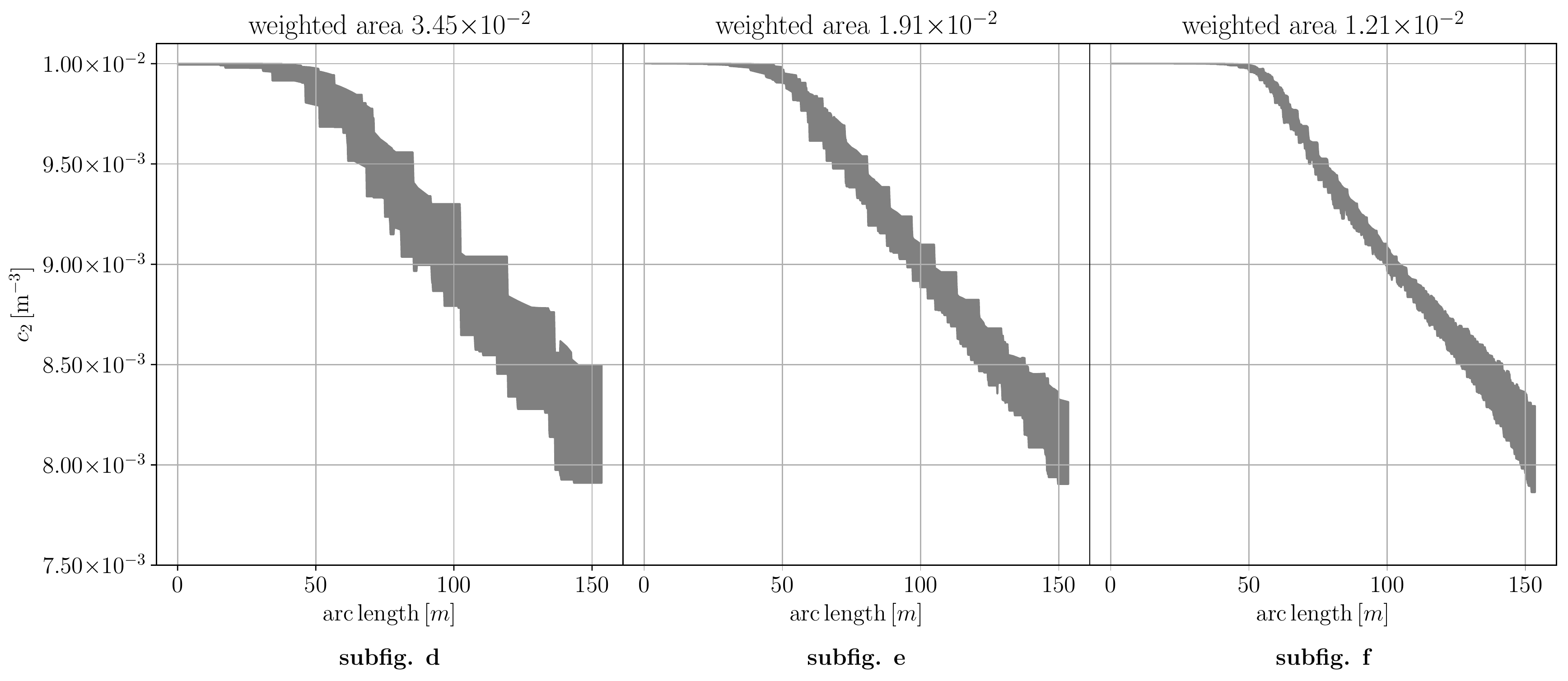}
    \caption{Case 1  of Subsection \ref{subsec:single_fracture}. On the top, concentration $c_2$ within the fracture, at the final simulation time, along the line $\left( \SI{0}{\metre}, \SI{100}{\metre}, \SI{80}{\metre} \right)$-$\left( \SI{100}{\metre}, \SI{0}{\metre}, \SI{20}{\metre} \right)$ for three refinements (coarse to fine). On the bottom, area between the 10th and 90th percentiles for three refinements (coarse to fine) and data. Results of Subsection \ref{case1:c}.}
    \label{fig:case1_pol_c_fracture}
\end{figure}

\paragraph{Fracture Concentration Over Line}\label{case1:c}
Figure \ref{fig:case1_pol_c_fracture} shows the concentration $c_2$ within the fracture at the final simulation time along the line $\left( \SI{0}{\metre}, \SI{100}{\metre}, \SI{80}{\metre} \right)$-$\left( \SI{100}{\metre}, \SI{0}{\metre}, \SI{20}{\metre} \right)$.

Again, almost all methods appear to converge with increasing refinement. \Ncu and \UibTpfa exhibit the largest deviations over all refinement levels. Close to the outlet boundary, \Ethz yields rather different values than the rest of the methods, but its tendency to approach the other methods with refinement can be observed clearly. Additionally, \Inm still shows considerably different results on the right boundary for the highest refinement level. Looking at the bottom row of Figure \ref{fig:case1_pol_c_fracture}, the convergence behavior of the spread is better than that of the matrix concentration reported in Subsection \ref{case1:b}, yet worse than for the matrix hydraulic head in Subsection \ref{case1:a}.

\begin{figure}[hbtp]
	\centering
    \includegraphics[width=0.95\textwidth]{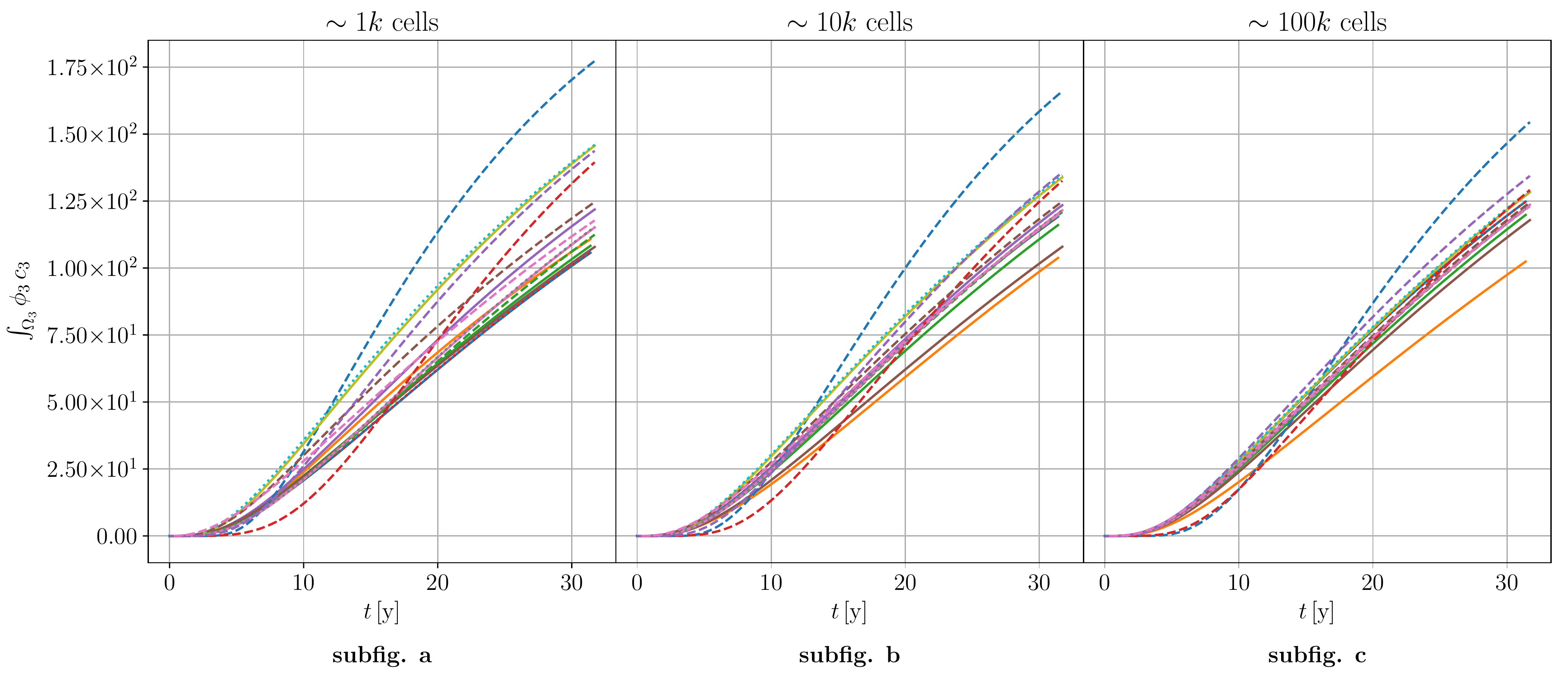}
    \caption{Case 1  of Subsection \ref{subsec:single_fracture}. integrated matrix concentration $\int_{\globalDomain_{3,3}} \porosity_3 c_3 \, \mathrm{d}x$ for three refinements (coarse to fine). Results of Subsection \ref{case1:d}.}
    \label{fig:case1_pot_c_matrix}
\end{figure}

\paragraph{Integrated Matrix Concentration Over Time}\label{case1:d}
Unlike the first three plots in \ref{case1:a}-\ref{case1:c}, Figure \ref{fig:case1_pot_c_matrix} illustrates an integrated quantity over time, namely, the integrated matrix concentration $\int_{\globalDomain_{3,3}} \porosity_3 c_3 \, \mathrm{d}x$.
Correspondingly, all curves appear much smoother than above. Over the three refinement levels, most methods again exhibit decreasing differences between each other. Remarkably, the  \UibTpfa shows a pronounced underestimation that increases over time. This can be explained by the inconsistency of the employed two-point flux approximation on the tetrahedral grids. Additionally, the \Ncu and \Ethz again exhibit larger differences.

\begin{figure}[hbtp]
	\centering
    \includegraphics[width=0.95\textwidth]{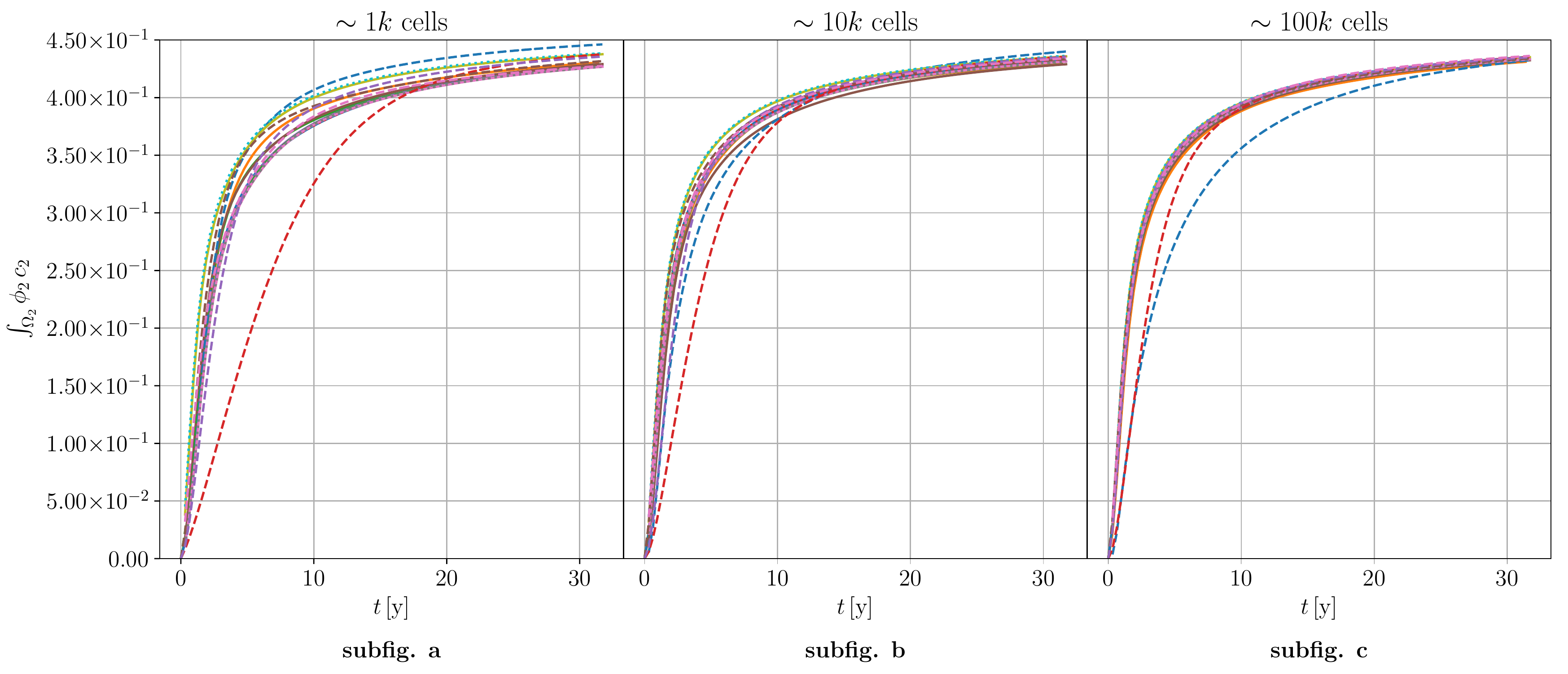}
    \caption{Case 1  of Subsection \ref{subsec:single_fracture}. Integrated fracture concentration $\int_{\globalDomain_2} \aperture_2 \porosity_2 c_2 \, \mathrm{d}x$ over time for three refinements (coarse to fine). Results of Subsection \ref{case1:e}.}
    \label{fig:case1_pot_c_fracture}
\end{figure}

\paragraph{Integrated Fracture Concentration Over Time}\label{case1:e}
Analogously, the integrated fracture concentration $\int_{\globalDomain_2} \aperture_2 \porosity_2 c_2 \, \mathrm{d}x$  for each time-step is visualized in Figure \ref{fig:case1_pot_c_fracture}.
The behavior of the curves is generally different from that reported in Subsection \ref{case1:d}, as the fracture fills up completely before the final simulation time. Here, the \UibTpfa is in line with the other methods whereas the \Ncu and \Ethz both deviate from the majority.

\begin{figure}[hbtp]
	\centering
    \includegraphics[width=0.95\textwidth]{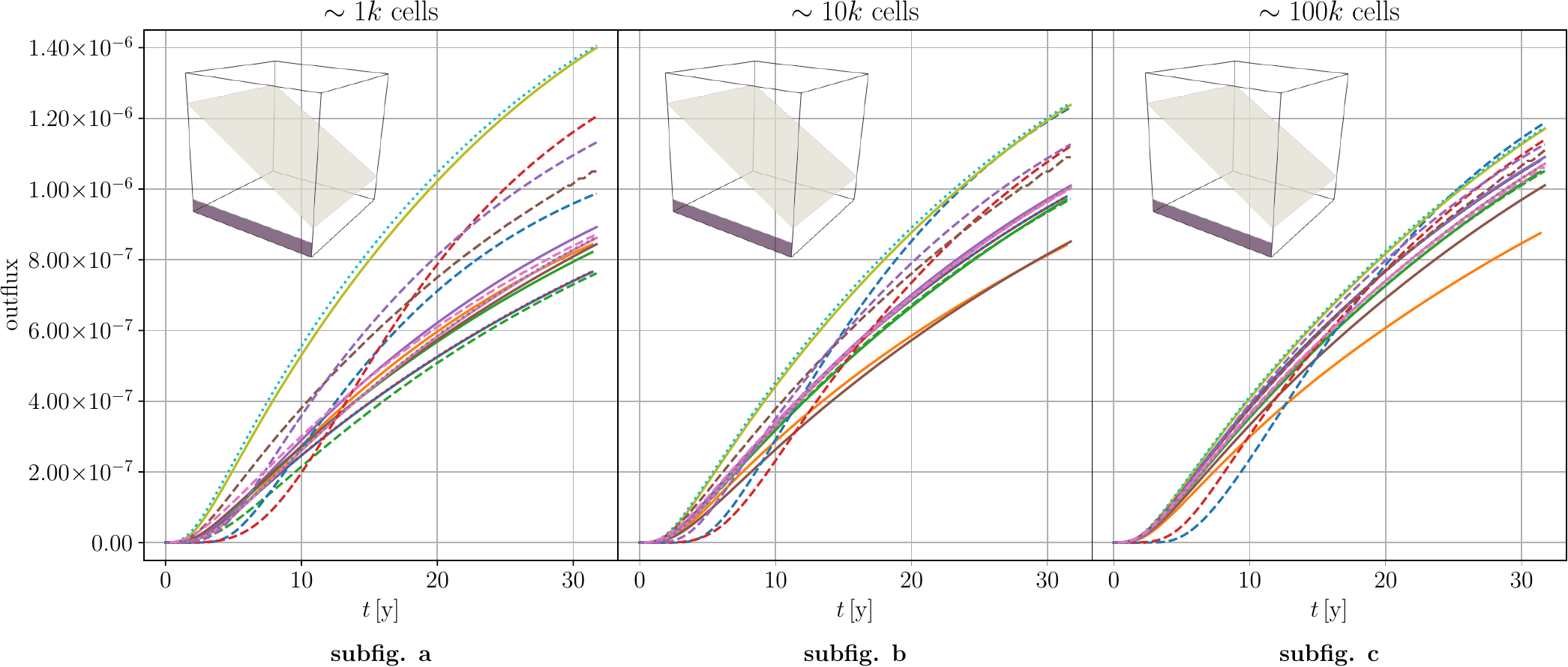}
    \caption{Case 1  of Subsection \ref{subsec:single_fracture}. Integrated flux of $c$ across the outlet boundary over time for three refinements (coarse to fine). Results of Subsection \ref{case1:f}.}
    \label{fig:case1_pot_outflux}
\end{figure}

\paragraph{Concentration Flux Across the Outlet Over Time}\label{case1:f}
Finally, Figure \ref{fig:case1_pot_outflux} depicts the integrated concentration flux across the outlet boundary over time.
Compared to the results in Subsection \ref{case1:e}, the agreement between the methods appears to be poorer. In particular, the two-point flux approximation of the \UibTpfa results in an underestimation similar to that reported in \ref{case1:d}. Again, \Ethz and \Ncu yield considerably different results at all refinement levels.

\paragraph{Computational Cost}\label{case1:g}
Indicators for the computational costs associated with the different methods are presented in Table \ref{tab:case1_comp_cost}. Most methods satisfy the prescribed numbers of elements. The most notable exception is given by the \Ncu, where six to ten times as many tetrahedral elements have been employed, to compensate for the fact that the degrees of freedom are associated with the vertices. The number of vertices are in line with the prescribed cell numbers. The relations of the number of degrees of freedom to the number of cells vary considerably between the different schemes, reflecting the characteristics from Table \ref{tab:scheme_properties_numerical}. The lowest such numbers are for the purely head- and vertex-based schemes on tetrahedrons for the \Ncu and \Dtu, while the highest ones result from the schemes that have head and velocity values as degrees of freedom. Additionally, the ratios of the number of nonzero entries to the number of degrees of freedom exhibit a large variability, ranging from approximately 5 (TPFA on tetrahedrons) to 30 (MPFA schemes with only head degrees of freedom).



\subsection{Case 2: Regular Fracture Network}\label{subsec:regular}

\noindent
\textbf{Benchmark case designers:} A. Fumagalli and I. Stefansson\\
\textbf{Benchmark case coordinators:} W. Boon and D. Gläser
\begin{figure}[hbtp]
	\centering
    \includegraphics[width=0.99\textwidth]{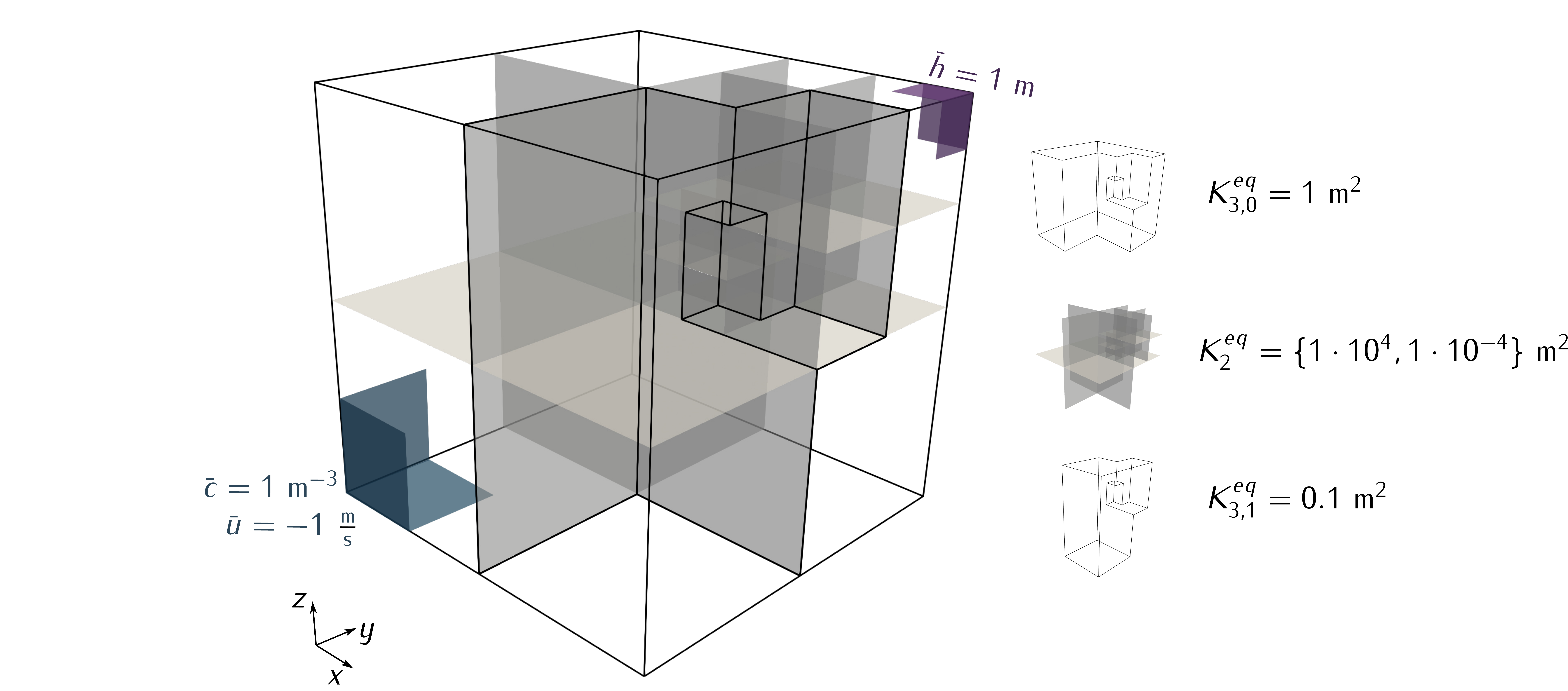}%
    \caption{Representation of the domain ($\Omega_3 = \left(0,1\right)^3$) and the fractures for 
             Case 2 of Subsection \ref{subsec:regular}. The inlet and outlet
             boundaries are colored in blue and purple, respectively, and on
             the right side, the permeability distributions among $\globalDomain_{3}$
             and $\globalDomain_{2}$ are illustrated.}
    \label{fig:domain_regular}
\end{figure}

\subsubsection{Description}
The second benchmark is a three-dimensional analog of the two-dimensional test case 4.1 from the 
benchmark study \cite{Flemisch:2018:BSF}. The domain is given by the unit cube $\Omega = \left(\SI{0}{\metre}, \SI{1}{\metre} \right)^3$ and contains 9 regularly oriented fractures, as illustrated in Figure \ref{fig:domain_regular}. 
The boundary $\partial \globalDomain$ is decomposed into three parts, each corresponding to a chosen boundary condition (see Figure \ref{fig:domain_regular}). First, $\partial \globalDomain_{h} = \{ (x, y, z) \in \partial \Omega: x, y, z > \SI{0.875}{\metre}\}$ is the part of the boundary on which we impose $\overline{h} = \SI{1}{\metre}$. Second, we set a flux boundary condition on $\partial \globalDomain_{in} = \{ (x, y, z) \in \partial \Omega: x, y, z < \SI{0.25}{\metre}\}$ by imposing $\overline{u} = \SI{-1}{\metre\per\second}$.
On the remainder of the boundary of $\globalDomain$, we impose no-flow conditions. 

Two variants of the test case are considered: {Case 2.1} has highly conductive fractures and {Case 2.2} has blocking fractures. In both cases, different hydraulic conductivities are prescribed in the following matrix subregions:
\begin{align*}
    \globalDomain_{3,0} =&\ \globalDomain_{3} \setminus \globalDomain_{3,1} \\
    \globalDomain_{3,1} =&\ \{(x, y, z) \in \globalDomain_{3}: x > \SI{0.5}{\metre} \cap y < \SI{0.5}{\metre}\} \\
    &\cup \{(x, y, z) \in \globalDomain_{3}: x > \SI{0.75}{\metre} \cap \SI{0.5}{\metre} < y < \SI{0.75}{\metre} \cap z > \SI{0.5}{\metre}\}\\
    &\cup \{(x, y, z) \in \globalDomain_{3}: \SI{0.625}{\metre} < x < \SI{0.75}{\metre} \cap \SI{0.5}{\metre} < y < \SI{0.625}{\metre} \cap \SI{0.5}{\metre} < z < \SI{0.75}{\metre}\}.
\end{align*}
For an illustration of these regions, we refer to the right part of Figure \ref{fig:domain_regular}. A complete overview of the parameters used in this test case is given in Table \ref{tab:case2_params}.

Finally, for the transport problem, we impose unitary 
concentration at the inflow boundary $\partial \globalDomain_{in}$.

\begin{table}[hbt]
\centering
\begin{tabular}{|l|ll|ll|}\hline
                        & \multicolumn{2}{l|}{{Case 2.1}} &  \multicolumn{2}{l|}{{Case 2.2}}      \\\hline
Matrix hydraulic conductivity $K_{3}|_{\Omega_{3, 0}}$        & $\bm{I}$& \si{\meter\per\second} & $\bm{I}$&\si{\meter\per\second}\\
Matrix hydraulic conductivity $K_{3}|_{\Omega_{3, 1}}$         & \num{1e-1}$\bm{I}$&\si{\meter\per\second}& \num{1e-1}$\bm{I}$&\si{\meter\per\second} \\
Fracture effective tangential hydraulic conductivity $K_2$ & $\bm{I}$ & \si{\meter\squared\per\second} & \num{1e-8}$\bm{I}$ & \si{\meter\squared\per\second}\\
Fracture effective normal hydraulic conductivity $\kappa_2$ & \num{2e8} & \si{\per\second} & $2$& \si{\per\second}\\
Intersection effective tangential hydraulic conductivity $K_1$ & \num{1e-4} &\si{\meter\cubed\per\second} & \num{1e-12} & \si{\meter\cubed\per\second}\\
Intersection effective normal hydraulic conductivity $\kappa_1$ & \num{2e4} &\si{\meter\per\second} & \num{2e-4} & \si{\meter\per\second}\\
Intersection effective normal hydraulic conductivity $\kappa_0$ & \num{2} &\si{\meter\squared\per\second} & \num{2e-8}& \si{\meter\squared\per\second}\\
Matrix porosity $\phi_{3}$ & \num{1e-1} & & \num{1e-1} & \\
Fracture porosity $\phi_2$ & \num{9e-1} & & \num{1e-2}& \\
Intersection porosity $\phi_1$ & \num{9e-1} & & \num{1e-2}& \\
Fracture cross-sectional length $\epsilon_2$  & \num{1e-4} &\si{\meter}& \num{1e-4} &\si{\meter}\\
Intersection cross-sectional area $\epsilon_1$  & \num{1e-8} &\si{\meter\squared} & \num{1e-8} &\si{\meter\squared}\\
Intersection cross-sectional volume $\epsilon_0$ & \num{1e-12} & \si{\meter\cubed} & \num{1e-12} &\si{\meter\cubed}\\
Total simulation time & \num{2.5e-1} & \si{\second} & &\\
Time-step $\Delta t$  & \num{2.5e-3} & \si{\second}  & &\\ \hline
\end{tabular}
\caption{Parameters used in Case 2 of Subsection \ref{subsec:regular}.}
\label{tab:case2_params}
\end{table}

\subsubsection{Results}
The results were collected for a sequence of 3 simulations by discretizing the 3d domain using approximately \num{500}, \num{4}k, and \num{32}k cells. The number of cells and degrees of freedom used by the participating methods are reported in Table \ref{tab:case2_comp_cost}. In the following, we discuss the results on the basis of line profiles of the hydraulic head in the 3d matrix as well as plots of the average concentrations within specified subregions of the 3d matrix. 

\begin{figure}[htbp]
	\centering
	\includegraphics[width=0.95\textwidth]{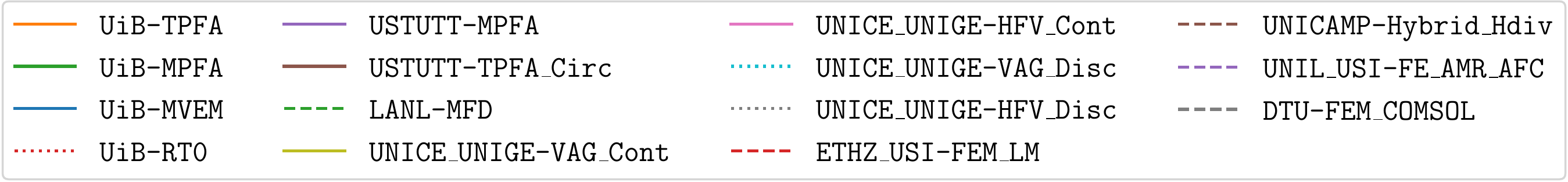}\\
	\includegraphics[width=0.95\textwidth]{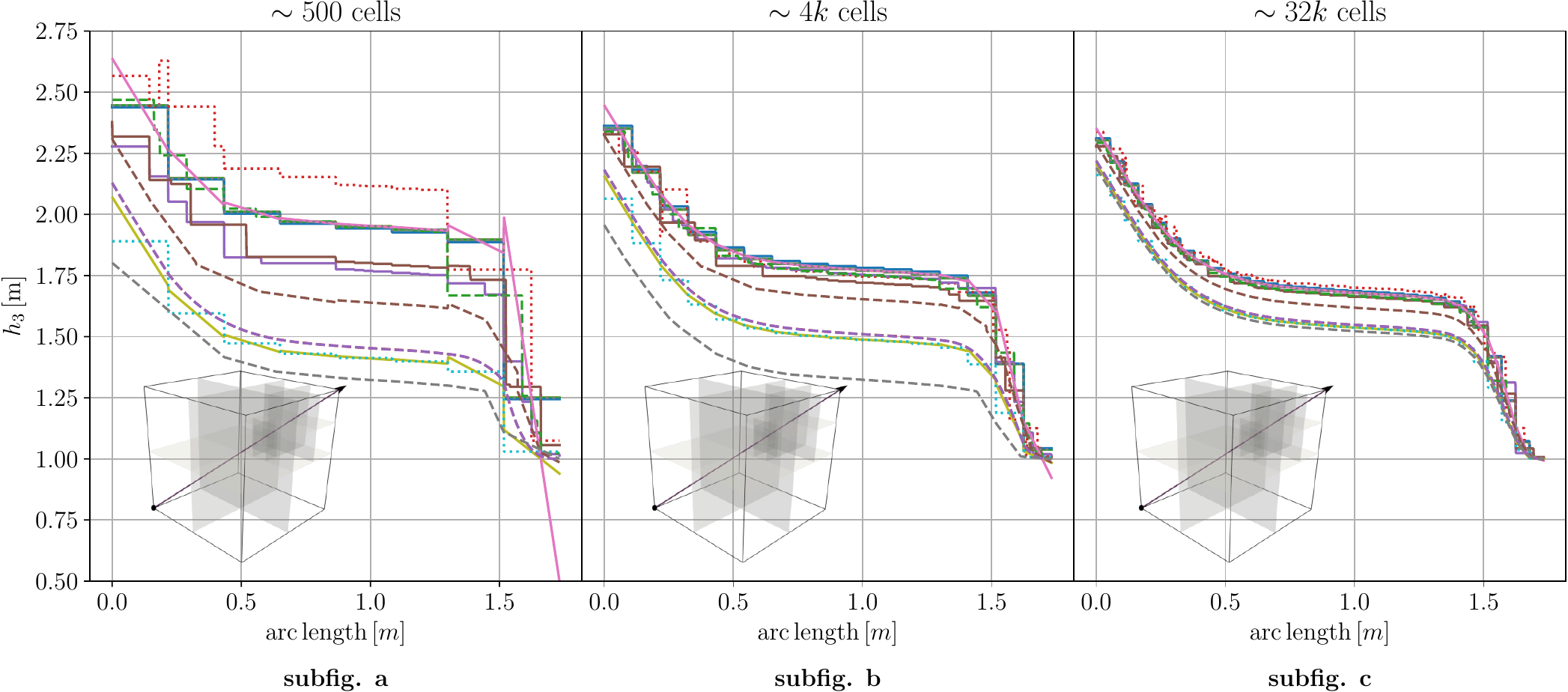}\\
    \includegraphics[width=0.95\textwidth]{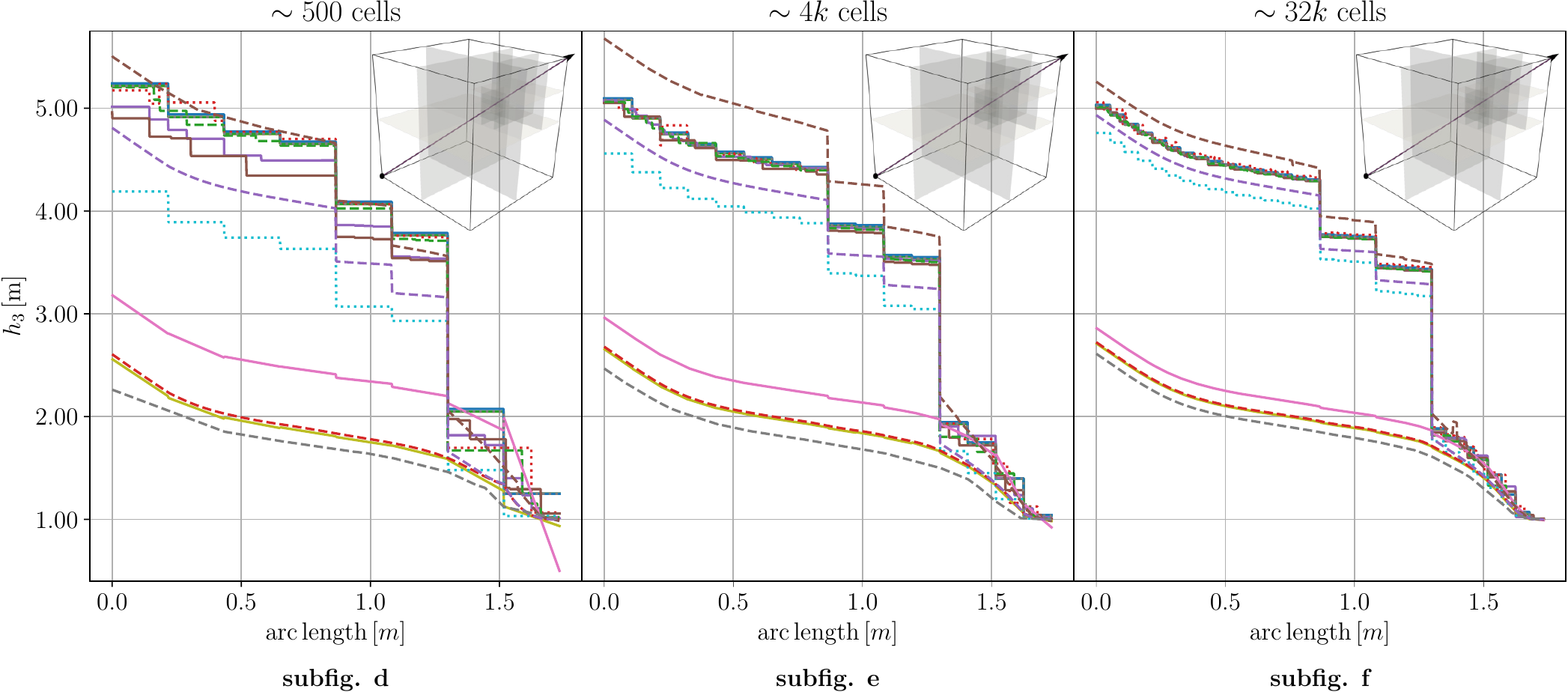}
    \caption{Case 2 of Subsection \ref{subsec:regular}. Plots of the hydraulic head $\head_3$ along the line $(0, 0, 0)$ - $(1, 1, 1)$ for the different
                     refinement levels (grid refinement increases from left to right) for the case of
                     conductive fractures ({Case 2.1}, upper row) and blocking fractures ({Case 2.2}, lower row). Results of Subsection \ref{case2:a}.}
    \label{fig:case2_pol}
\end{figure}

\paragraph{Hydraulic Head Over Line}\label{case2:a}
Figure \ref{fig:case2_pol} shows the hydraulic head $\head_3$ plotted along the diagonal line segment (0, 0, 0)-(1, 1, 1) for all grid refinements and for both {Case 2.1} and {Case 2.2}. In the case of conductive fractures
the spread decreases significantly upon grid refinement, although some noticeable differences still prevail for the finest grid.

In the case of blocking fractures, the highest discrepancies are shown by the schemes that assume continuity of the hydraulic head across the fractures. As expected, these methods cannot capture the jump in the hydraulic head present in this test case. On the other hand, the remaining schemes seem to approach the same solution. We observe that the \UniceVagD and the \Unil produce slightly lower and the \Unicamp scheme slightly higher hydraulic heads, but the deviations tend to diminish with increasing grid refinement.

The \UniceVagC and \UniceVagD methods incorporate Dirichlet boundary conditions on the
vertices rather than on faces.
This may explain, in part, the deviations in hydraulic head observed on coarse meshes for these methods. 
As expected, these differences decrease with mesh refinement.
For the \Unil method, the differences might come from the
representation of the fractures, which have the same spatial dimension as the background matrix. In particular, each fracture consists of a layer of elements that is refined at least twice by using adaptive mesh refinement.

\paragraph{Mean Matrix Concentration Over Time}\label{case2:b}
The second comparison in {Case 2} concerns the solution of the transport equation over time.  
These solutions are computed only on the second level of mesh refinement, i.e.,\ using approximately \num{4000} cells. For the simulation of the transport model, the upwind scheme is employed for all methods except \Unil and \Ethz, which employ a finite element discretization with an algebraic flux correction \cite{kuzmin2012flux}.

\begin{figure}[htbp]
	\centering
	\includegraphics[width=0.95\textwidth]{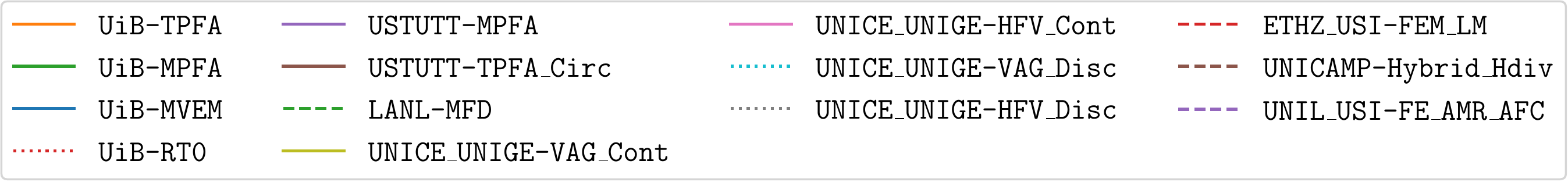}\\
	\includegraphics[width=0.95\textwidth]{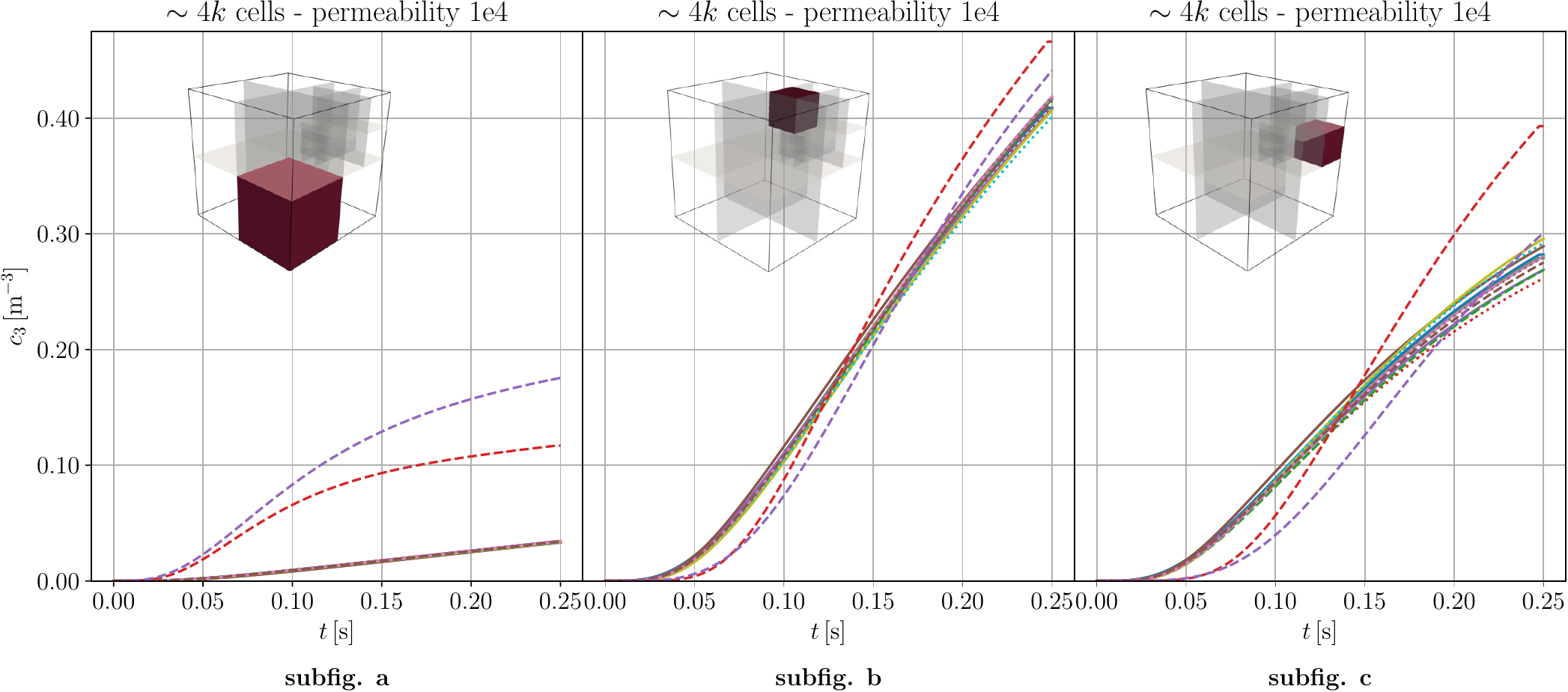}\\
    \includegraphics[width=0.95\textwidth]{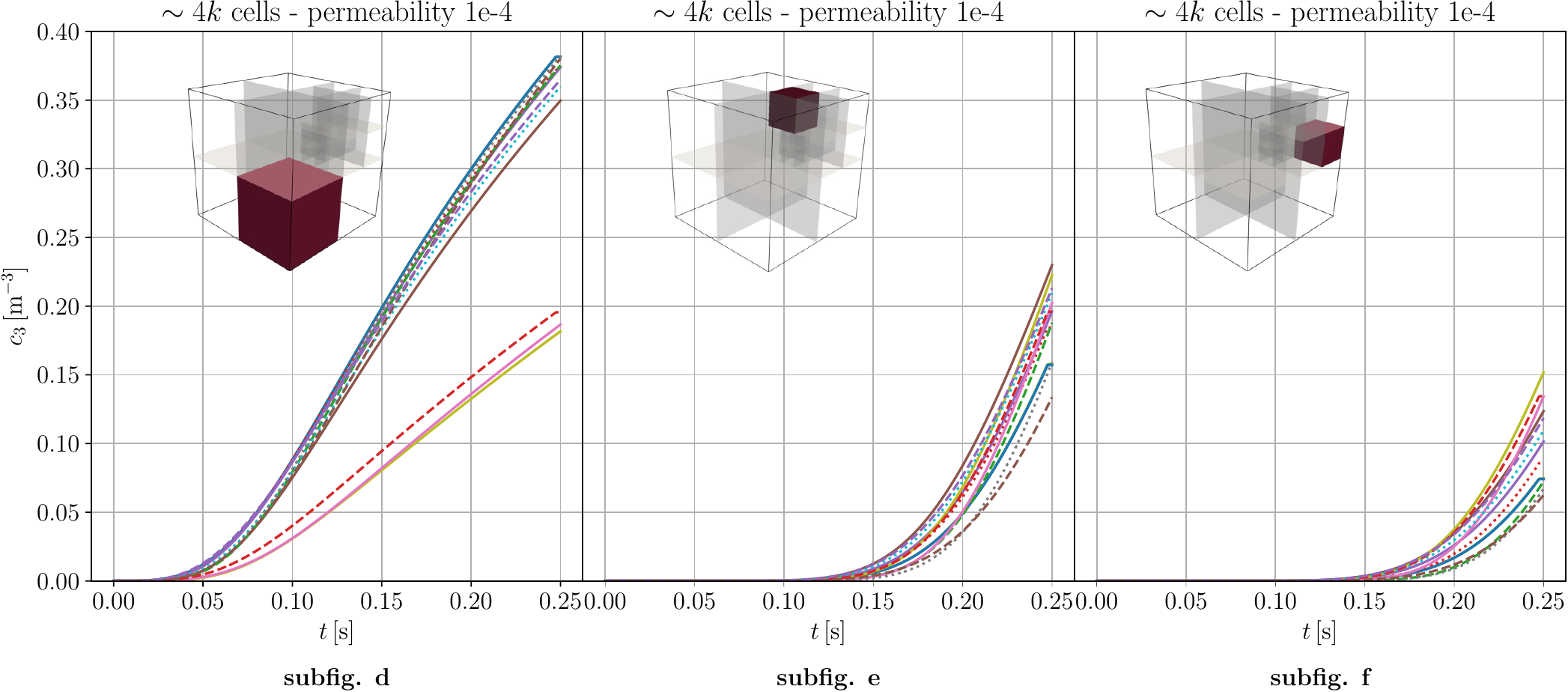}
    \caption{Case 2  of Subsection \ref{subsec:regular}. On the top, temporal evolution of the average tracer concentration in matrix regions $1$, $10$ and $11$ (from left to right) for the case of conductive fractures ({Case 2.1}). On the bottom, temporal evolution of the average tracer concentration in the matrix regions $1$, $10$ and $11$ (from left to right) for the case of blocking fractures ({Case 2.2}).  Results of Subsection \ref{case2:b}.}
    \label{fig:case2_cot_cond}
\end{figure}

The top of Figure \ref{fig:case2_cot_cond} depicts the temporal evolution of the mean tracer concentrations in three matrix regions for the case of highly conductive fractures. These regions were selected to form a representative illustration of the spread between the schemes.
It can be seen that the majority of the schemes produce rather low concentrations in the first region, on the order of $\SI{2.5}{\percent}$ at the final simulation time. In contrast, the \Ethz and the \Unil schemes produce significantly higher concentrations with values above $\SI{10}{\percent}$ at the end of the simulation. In general, the temporal evolution of the concentrations in these three regions agrees very well among the majority of participating schemes, while the \Ethz and the \Unil schemes show significant deviations.
These might be related to the flow discretization methods, but could also be affected by the different discretization that is employed for the transport discretization related to these methods, and, for \Unil, also the underlying equidimensional model.

For the case of blocking fractures, the concentrations in the same matrix regions are illustrated in the bottom row of Figure \ref{fig:case2_cot_cond}. In general, a larger spread of the computed concentrations can be observed. For the first region, the schemes that assume continuity of the hydraulic head produce significantly lower concentrations, while the remaining schemes produce solutions that agree rather well. However, for the second and third regions, the concentrations at the final simulation time show a widespread among all participating schemes. 

As a general trend, it can be observed that the differences in computed concentrations increase with time. Additionally, differences increase with the regions' distance from the inflow boundary. As expected, for the case of conductive fractures, the differences are smaller than in the case of blocking fractures.   
%



\subsection{Case 3: Network with Small Features} \label{subsec:small_features}

\noindent
\textbf{Benchmark case designers:} E. Keilegavlen and I. Stefansson\\
\textbf{Benchmark case coordinator:} I. Stefansson and A. Fumagalli
\subsubsection{Description}
\begin{figure}[hbtp]
\centering
\includegraphics[width=0.9\textwidth]{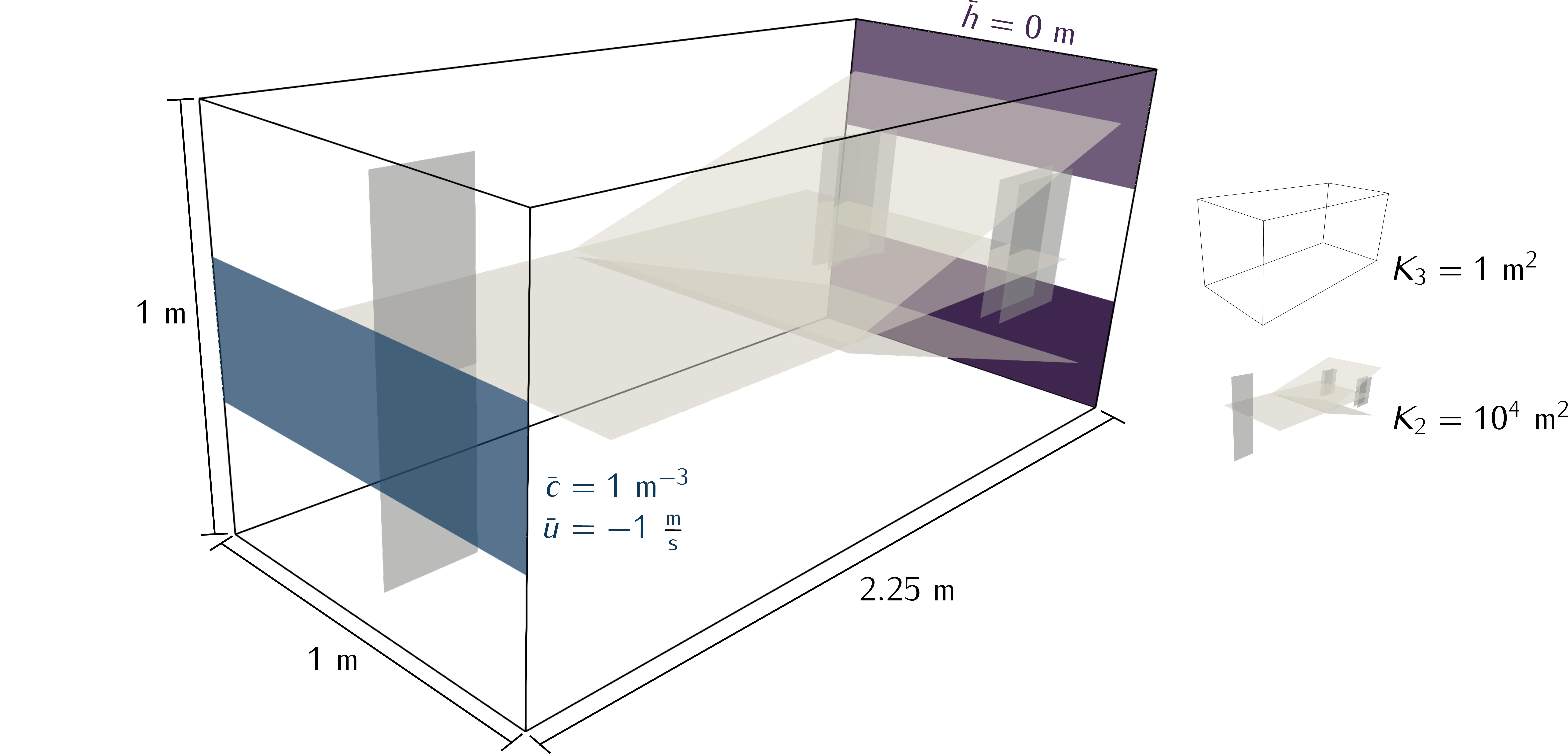}\hfill%
\caption{Representation of the fractures and the outline of the domain
    	for Case 3 of Subsection \ref{subsec:small_features}.}
\label{fig:geom_small_features}
\end{figure}

This test case is designed to probe accuracy in the presence of small geometric features, which may cause trouble for conforming meshing strategies.
The domain is the box $\Omega = (\SI{0}{\metre}, \SI{1}{\metre}) \times (\SI{0}{\metre}, \SI{2.25}{\metre}) \times (\SI{0}{\metre}, \SI{1}{\metre})$, containing eight fractures (see Figure \ref{fig:geom_small_features}).

We define the inlet and outlet boundaries as follows:
\begin{align*}
\partial\globalDomain_N &= \partial\Omega \setminus (\partial\globalDomain_{in} \cup \partial\globalDomain_{out}) \\
\partial\globalDomain_{in} &= (\SI{0}{\metre}, \SI{1}{\metre}) \times \{\SI{0}{\metre}\} \times (\SI{1/3}{\metre}, \SI{2/3}{\metre})\\
\partial\globalDomain_{out} &= \partial\globalDomain_{out, 0} \cup \partial\globalDomain_{out, 1}\\
\partial\globalDomain_{out, 0} &= (\SI{0}{\metre}, \SI{1}{\metre}) \times \{\SI{2.25}{\metre}\} \times (\SI{0}{\metre}, \SI{1/3}{\metre}) \\
\partial\globalDomain_{out, 1} &= (\SI{0}{\metre}, \SI{1}{\metre}) \times \{\SI{2.25}{\metre}\} \times (\SI{2/3}{\metre}, \SI{1}{\metre})
\end{align*}
The boundary conditions for flow are homogeneous Dirichlet conditions on $\partial\globalDomain_{out}$, uniform unit inflow on $\partial\globalDomain_{in}$, so that $\int_{\partial\globalDomain_{in}} \vecu_3 \cdot\normal dS = -1/3$ \si{\metre\cubed\per\second}, and homogeneous Neumann conditions on $\partial\globalDomain_N$. For the transport problem, we consider
a homogeneous initial condition and as boundary condition a unit concentration at $\partial\globalDomain_{in}$. A complete overview of the parameters used in Case 3 is given in Table \ref{tab:case3_params}.

\begin{table}[hbt]
\centering
\begin{tabular}{|l|ll|}\hline
Matrix hydraulic conductivity $K_3$& $\bm{I}$ & \si{\metre\per\second} \\
Fracture effective tangential hydraulic conductivity $K_2$ & \num{1e2}$\bm{I}$ & \si{\metre\squared\per\second} \\
Fracture effective normal hydraulic conductivity $\kappa_2$ & \num{2e6} & \si{\per\second} \\
Intersection effective tangential hydraulic conductivity $K_1$ & $1$ &  \si{\metre\cubed\per\second} \\
Intersection effective normal hydraulic conductivity $\kappa_1$ & \num{2e4} & \si{\metre\per\second} \\
Matrix porosity $\phi_3$ & \num{2e-1} & \\
Fracture porosity $\phi_2$ & \num{2e-1} & \\
Intersection effective porosity $\phi_1$ & \num{2e-1} &  \\
Fracture cross-sectional length $\epsilon_2$ & \num{1e-2} & \si{\metre} \\
Intersection cross-sectional area $\epsilon_1$ & \num{1e-4} & \si{\metre\squared} \\
Total simulation time & \num{1e0} & \si{\second} \\
Time-step $\Delta t$  & \num{1e-2} & \si{\second} \\ \hline
\end{tabular}
\caption{Parameters used in Case 3 of Subsection \ref{subsec:small_features}.}
\label{tab:case3_params}
\end{table}

\subsubsection{Results} 
 Similar to the previous cases, we compare the methods on the basis of a) the hydraulic head of the matrix domain along two lines, b) the integrated fracture concentration over time, c) the fluxes out of the domain and d) computational cost.
Two different simulations with approximately \num{30}k and \num{150}k cells for the 3d domain were performed.
It was seen as infeasible to include one more level of refinement for all methods. However, refined versions of the \StuttMpfa with up to approximately \num{1e6} matrix cells were produced. 
At this stage, there were no noticeable differences between solutions on different grids, and the finest solution was included as a reference solution.

\begin{figure}[hbt]
	\centering
    \includegraphics[width=0.95\textwidth]{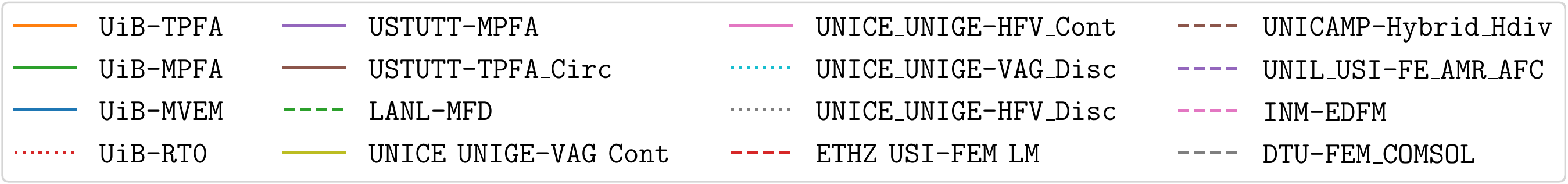}\\
	\includegraphics[width=0.95\textwidth]{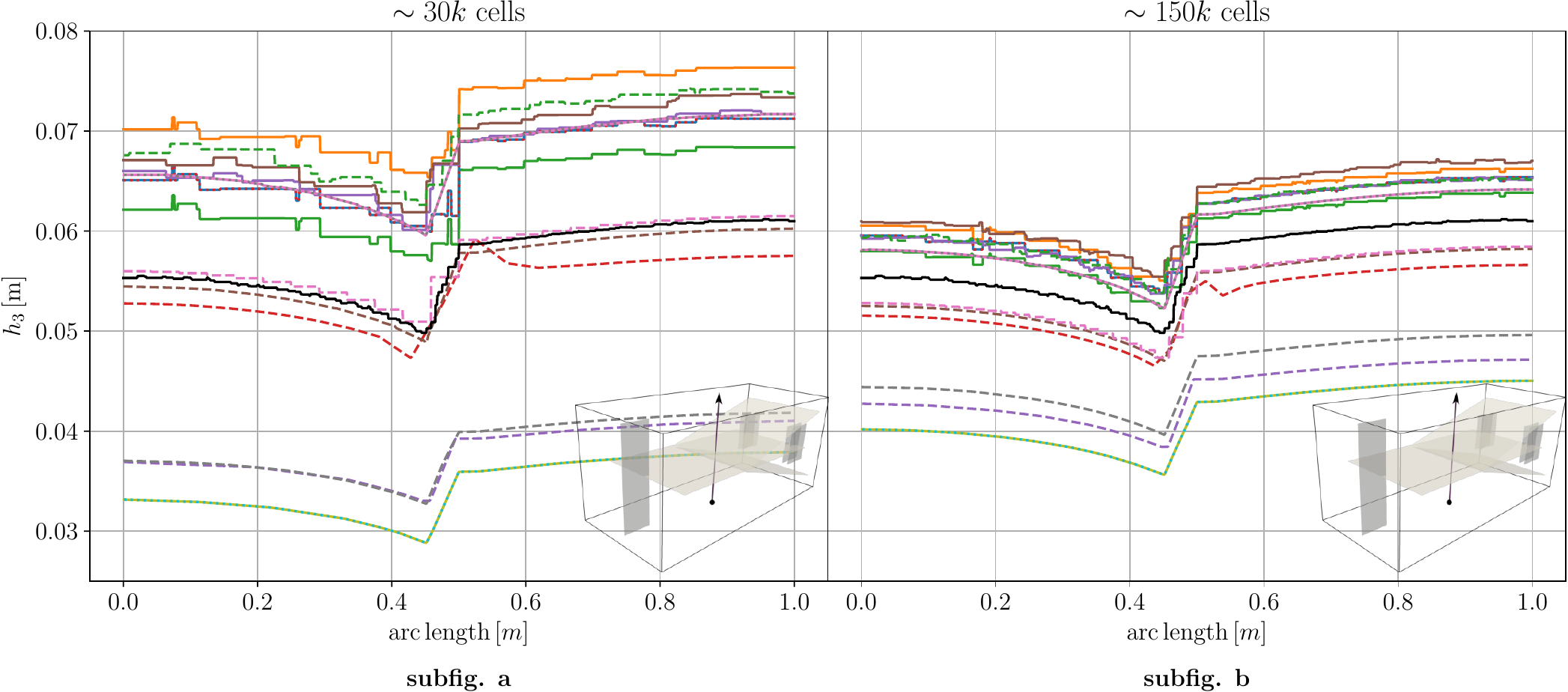}
    \caption{Case 3 of Subsection \ref{subsec:small_features}. Hydraulic head $\head_3$ in the matrix over the line $\left(\SI{0.5}{\metre}, \SI{1.1}{\metre}, \SI{0}{\metre} \right)$-$\left(\SI{0.5}{\metre}, \SI{1.1}{\metre}, \SI{1}{\metre} \right)$ for the coarse (left) and fine (right) grid. The solid black line shows the solution obtained with the \StuttMpfa scheme on a grid with approximately \num{1e6} matrix cells. Results of Subsection \ref{case3:a}.}%
    \label{fig:case3_pol_p_0_matrix}
\end{figure}

\paragraph{Hydraulic Head Over Line}\label{case3:a}
Figure \ref{fig:case3_pol_p_0_matrix} shows the profile of the hydraulic head $\head_3$ in the matrix along the line $\left(\SI{0.5}{\metre}, \SI{1.1}{\metre}, \SI{0}{\metre} \right)$-$\left(\SI{0.5}{\metre}, \SI{1.1}{\metre}, \SI{1}{\metre} \right)$.
This shows considerable differences between the methods for both refinement levels.
However, the agreement is better for the second refinement level, where most of the methods are within a relative hydraulic head range of approximately \SI{10}{\percent}. The \UniceVagD, \UniceVagC, \Dtu, and \Unil methods show the highest discrepancies in these plots, but the deviation from the reference solution decreases significantly with higher refinement.
The significant difference between the refinements may indicate that the small features of the fracture network geometry are not adequately resolved, at least not by the coarser grids. This is in line with the purpose of the test case.

\begin{figure}[hbtp]
	\centering
    \includegraphics[width=0.95\textwidth]{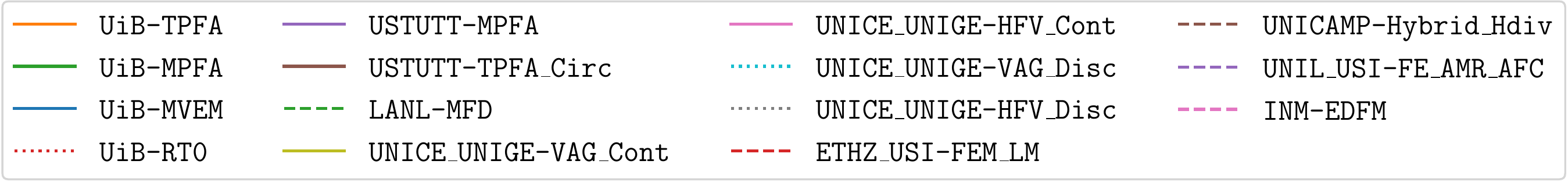}\\
    \includegraphics[width=0.95\textwidth]{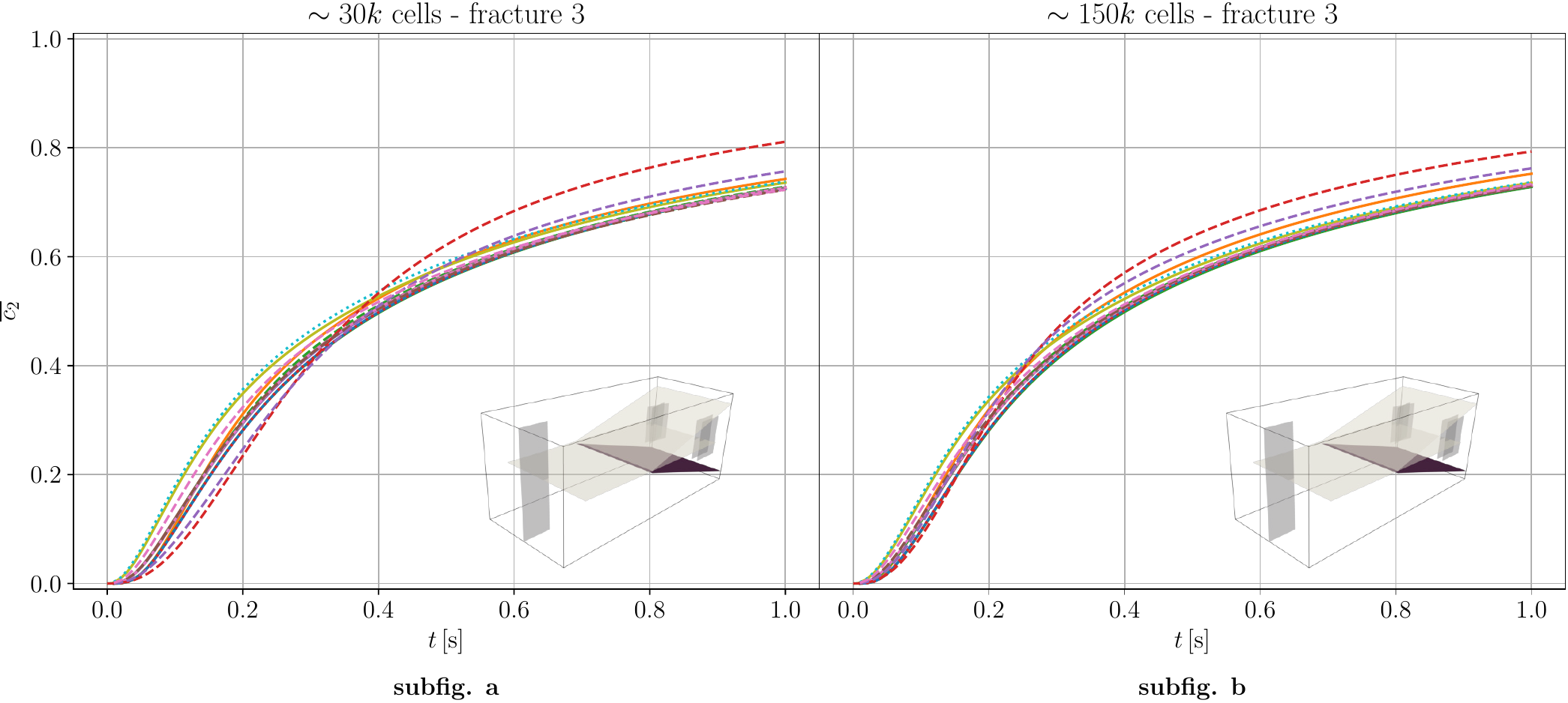}
    \caption{Case 3 of Subsection \ref{subsec:small_features}. Mean concentration within fracture number 3 throughout the simulation time for the coarse (left) and fine (right) grid. Results of Subsection \ref{case3:b}.}%
    \label{fig:case3_pot_fracture_3}
\end{figure}

\paragraph{Mean Fracture Concentration Over Time}\label{case3:b}
Data were reported for the integrated concentration $\overline{c_2} = \int_{\Omega_{2,i}}c_2 / |\Omega_{2,i}|$ on each fracture $i$ throughout the simulation. 
There is a general agreement between the methods, with the method of \Ethz showing some deviations for some of the fractures. As an example, Figure \ref{fig:case3_pot_fracture_3} shows the plots for both refinement levels for fracture number 3, demonstrating limited difference between the refinement levels.

\begin{figure}[hbtp]
	\centering
    \includegraphics[width=0.95\textwidth]{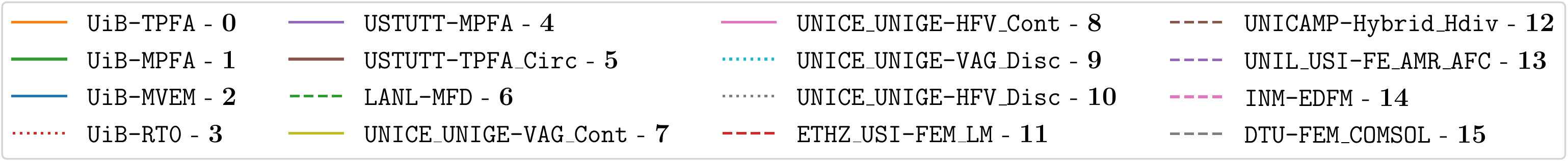}\\
    \includegraphics[width=0.45\textwidth]{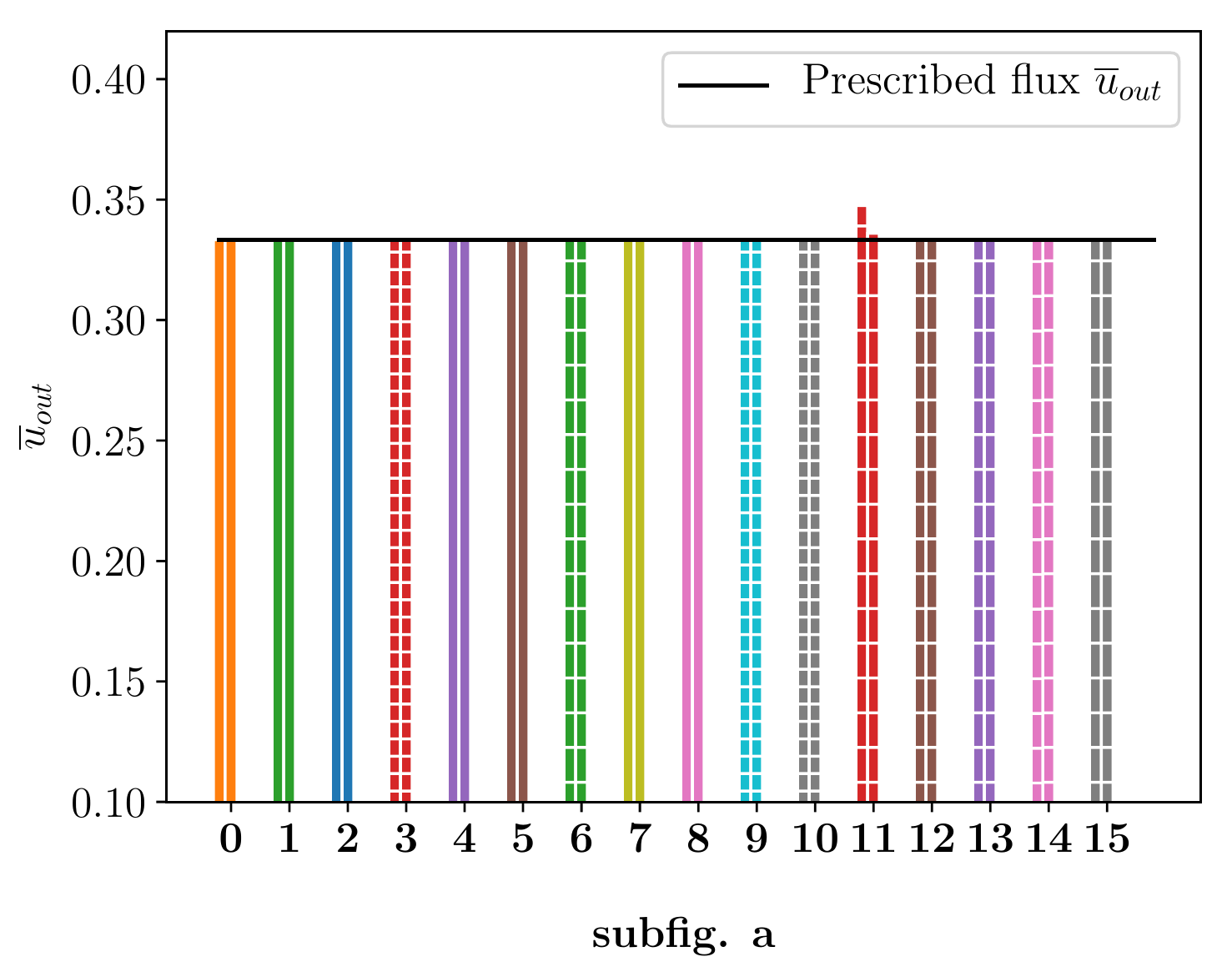}
    \includegraphics[width=0.45\textwidth]{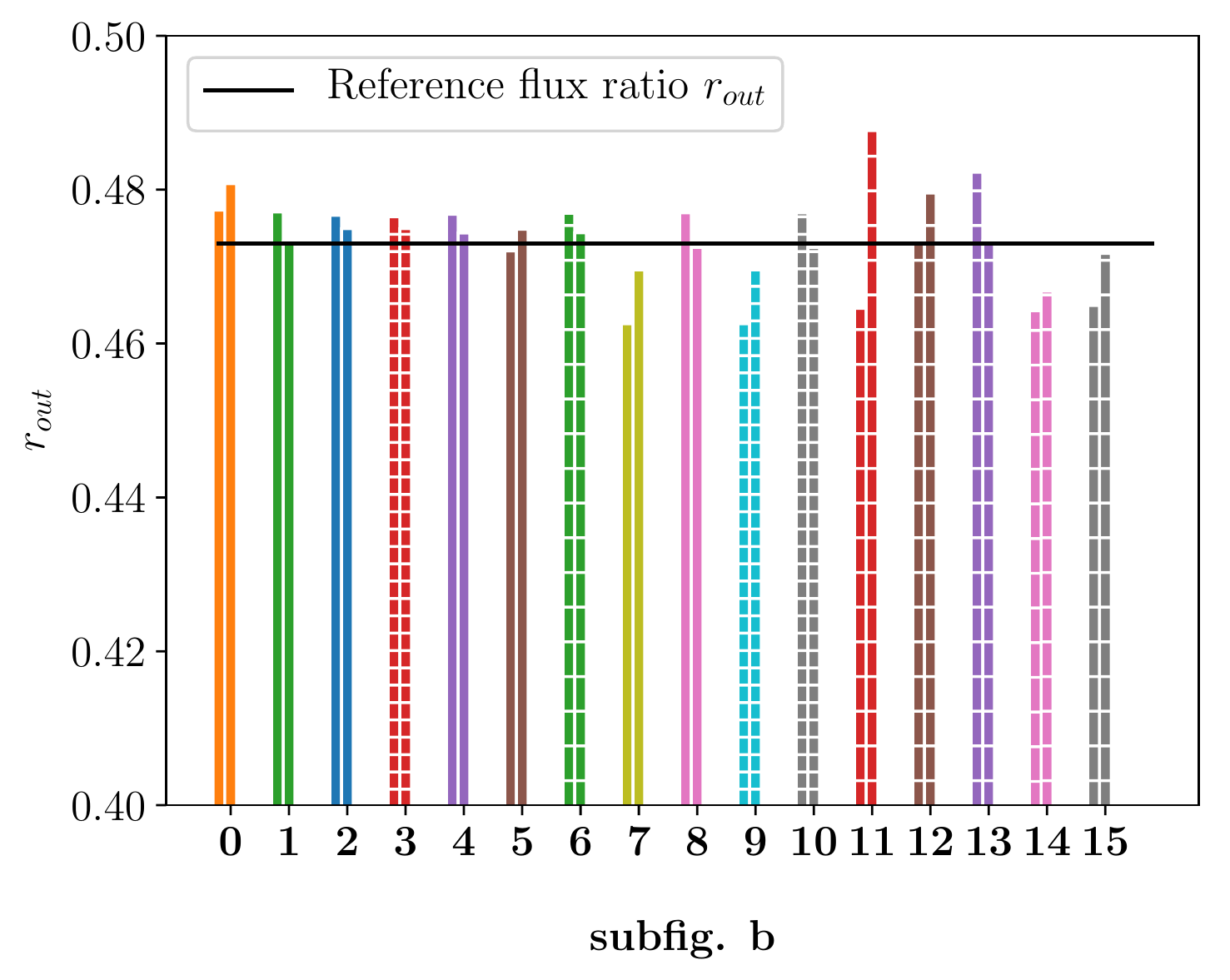}
    \caption{Case 3 of Subsection \ref{subsec:small_features}. Total outflux (left) and ratio exiting over $\partial\globalDomain_{out, 0}$ (right). The bar pairs correspond to the coarse and fine grid, while the reference solution is indicated by the horizontal line. Results of Subsection \ref{case3:c}.}%
    \label{fig:case3_boundary_data}
\end{figure}

\paragraph{Boundary Fluxes}\label{case3:c}
The total outflow $\overline{u}_{out}=\int_{\partial\globalDomain_{out}} \vecu_3 \cdot\normal dS$ and the proportion exiting over $\partial\globalDomain_{out, 0}$, i.e., $r_{out} =\int_{\partial\globalDomain_{out,0}} \vecu_3 \cdot\normal dS/\overline{u}_{out}$, are shown in Figure \ref{fig:case3_boundary_data}. When compared to the prescribed inflow of \SI{-1/3}{\metre\cubed\per\second}, the $\overline{u}_{out}$ values reveal a small lack of volume conservation for \Ethz, but the method improves for the finer grid. The ratio $r_{out}$ provides an indication of whether the flux fields agree. The ratios generally agree well with the refined \StuttMpfa, except for the \Ethz method, which does not approach the reference value for the finest grid.

\paragraph{Computational Cost}\label{case3:e}
Based on the data presented in Table \ref{tab:case3_comp_cost}, we note that the \Unil applies \num{68}k and \num{~203}k cells for the cases where \num{30}k and \num{150}k cells were prescribed, respectively. The rest of the methods are well within \SI{10}{\percent} of the prescribed values. As for the other test cases, there are significant variations in the number of degrees of freedom and nonzero matrix entries related to the design of the methods.


\subsection{Case 4: Field Case} \label{subsec:field_network}

\noindent
\textbf{Benchmark case designers:} E. Keilegavlen and A. Fumagalli\\
\textbf{Benchmark case coordinator:} E. Keilegavlen
\subsubsection{Description}
\begin{figure}[tb]
\centering
\includegraphics[width=.75\textwidth]{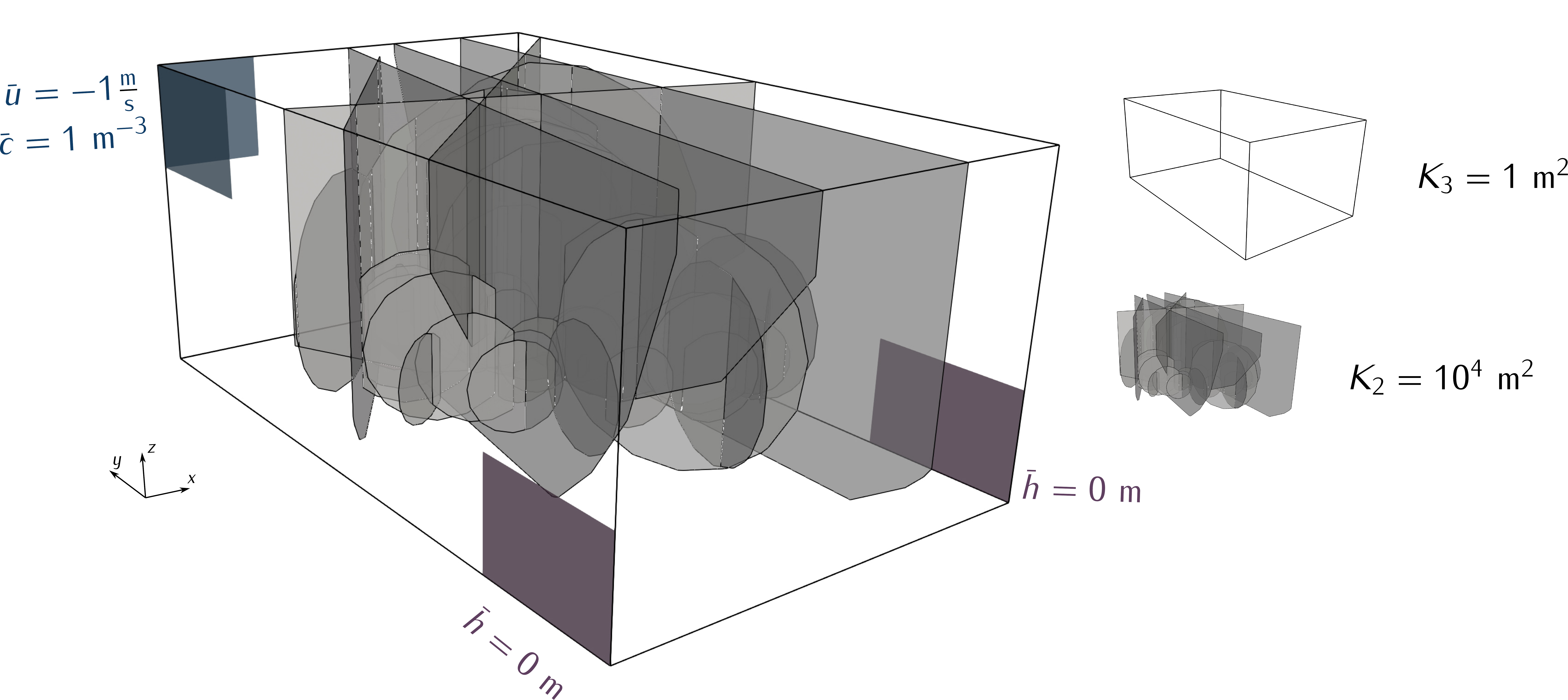}\hfill%
\caption{Case 4 of Subsection \ref{subsec:field_network}. Representation of the fractures and the outline of the domain. Inlet boundaries are shown in blue, outlets in purple.}
\label{fig:case4_network}
\end{figure}

The geometry of the fourth case is based on a postprocessed 
outcrop from the island of Algerøyna, outside Bergen, Norway, and is a subset of the fracture network presented in \cite{Fumagalli2017e}.
From the outcrop, 52 fractures were selected, extruded in the vertical direction and then cut by a bounding box. 
The resulting network has 106 fracture intersections, and multiple fractures intersect the domain boundary.
The simulation domain is the box $\Omega = (\SI{-500}{\metre}, \SI{350}{\metre}) \times (\SI{100}{\metre}, \SI{1500}{\metre}) \times (\SI{-100}{\metre}, \SI{500}{\metre})$. 
The fracture geometry is depicted in Figure \ref{fig:case4_network}.

The inlet and outlet boundaries are defined as follows:
\begin{align*}
\partial\globalDomain_N &= \partial\Omega \setminus (\partial\globalDomain_{in} \cup \partial\globalDomain_{out}), \\
\partial\globalDomain_{in} &= \partial\globalDomain_{in, 0} \cup \partial\globalDomain_{in, 1},\qquad
\partial\globalDomain_{out} = \partial\globalDomain_{out, 0} \cup \partial\globalDomain_{out, 1},\\
\partial\globalDomain_{in, 0} &= (\SI{-500}{\metre}, \SI{-200}{\metre}) \times \{\SI{1500}{\metre}\} \times (\SI{300}{\metre}, \SI{500}{\metre}), \\
\partial\globalDomain_{in, 1} &= \{\SI{-500}{\metre}\} \times (\SI{1200}{\metre}, \SI{1500}{\metre}) \times (\SI{300}{\metre}, \SI{500}{\metre}),\\
\partial\globalDomain_{out, 0} &= \{\SI{-500}{\metre}\} \times (\SI{100}{\metre}, \SI{400}{\metre}) \times (\SI{-100}{\metre}, \SI{100}{\metre}),\\
\partial\globalDomain_{out, 1} &= \{\SI{350}{\metre}\} \times (\SI{100}{\metre}, \SI{400}{\metre}) \times (\SI{-100}{\metre}, \SI{100}{\metre}).
\end{align*}
The boundary conditions for flow are homogeneous Dirichlet conditions on $\partial\globalDomain_{out}$, uniform unit inflow on $\partial\globalDomain_{in}$, so that $\int_{\partial\globalDomain_{in}} \vecu_3 \cdot\normal dS = \SI{-1.2e5}{\metre\cubed\per\second}$, and homogeneous Neumann conditions on $\partial\globalDomain_N$. For the transport problem, we consider a
homogeneous initial condition, with a unit concentration at $\partial\globalDomain_{in}$. The parameters for conductivity, porosity and aperture are given in Table \ref{tab:parameters_case4}, as is the total simulation time and time-step size.

\begin{table}[hbt]
\begin{center}
\begin{tabular}{|l|ll|}\hline
Matrix hydraulic conductivity $K_3$& $\bm{I}$ & \si{\metre\per\second} \\
Fracture effective tangential hydraulic conductivity $K_2$ & \num{1e2}$\bm{I}$ & \si{\metre\squared\per\second} \\
Fracture effective normal hydraulic conductivity $\kappa_2$ & \num{2e6} & \si{\per\second} \\
Intersection effective tangential hydraulic conductivity $K_1$ & $1$ &  \si{\metre\cubed\per\second} \\
Intersection effective normal hydraulic conductivity $\kappa_1$ & \num{2e4} & \si{\metre\per\second} \\
Matrix porosity $\phi_3$ & \num{2e-1} & \\
Fracture porosity $\phi_2$ & \num{2e-1} & \\
Intersection porosity $\phi_1$ & \num{2e-1} & \\
Fracture cross-sectional length $\epsilon_2$ & \num{1e-2} & \si{\metre} \\
Intersection cross-sectional area $\epsilon_1$ & \num{1e-4} & \si{\metre\squared} \\
Total simulation time & \num{5e3} & \si{\second} \\
Time-step $\Delta t$  & \num{5e1} & \si{\second} \\ \hline
\end{tabular}
\end{center}
\caption{Parameter used in Case 4  of Subsection \ref{subsec:field_network}.}
\label{tab:parameters_case4}
\end{table}

Because of the complex network geometry, grid refinement studies were considered infeasible and the benchmark specified the usage of a single grid. A Gmsh \cite{gmsh2009} configuration file, was provided to assist participants with geometry processing and meshing. The use of this predefined grid was optional, but the number of 3d cells should be approximately 260k.

\subsubsection{Results}
Results were reported for 14 schemes. 
The two methods that participated in Case 3, which is closest  in geometric complexity, but not in Case 4, are \Inm and \Unil.
The participating methods are compared in terms of a) hydraulic head of the matrix domain along two lines, b) time series of concentrations in selected fractures and c) computational cost.

\begin{figure}[tb]
\centering
\includegraphics[width=.95\textwidth]{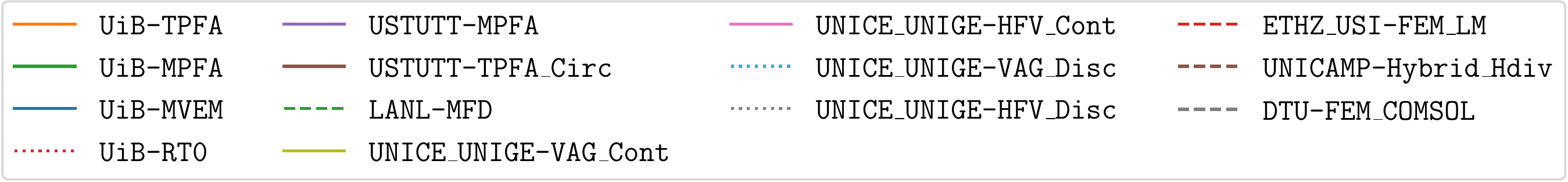}\\
\includegraphics[width=.45\textwidth]{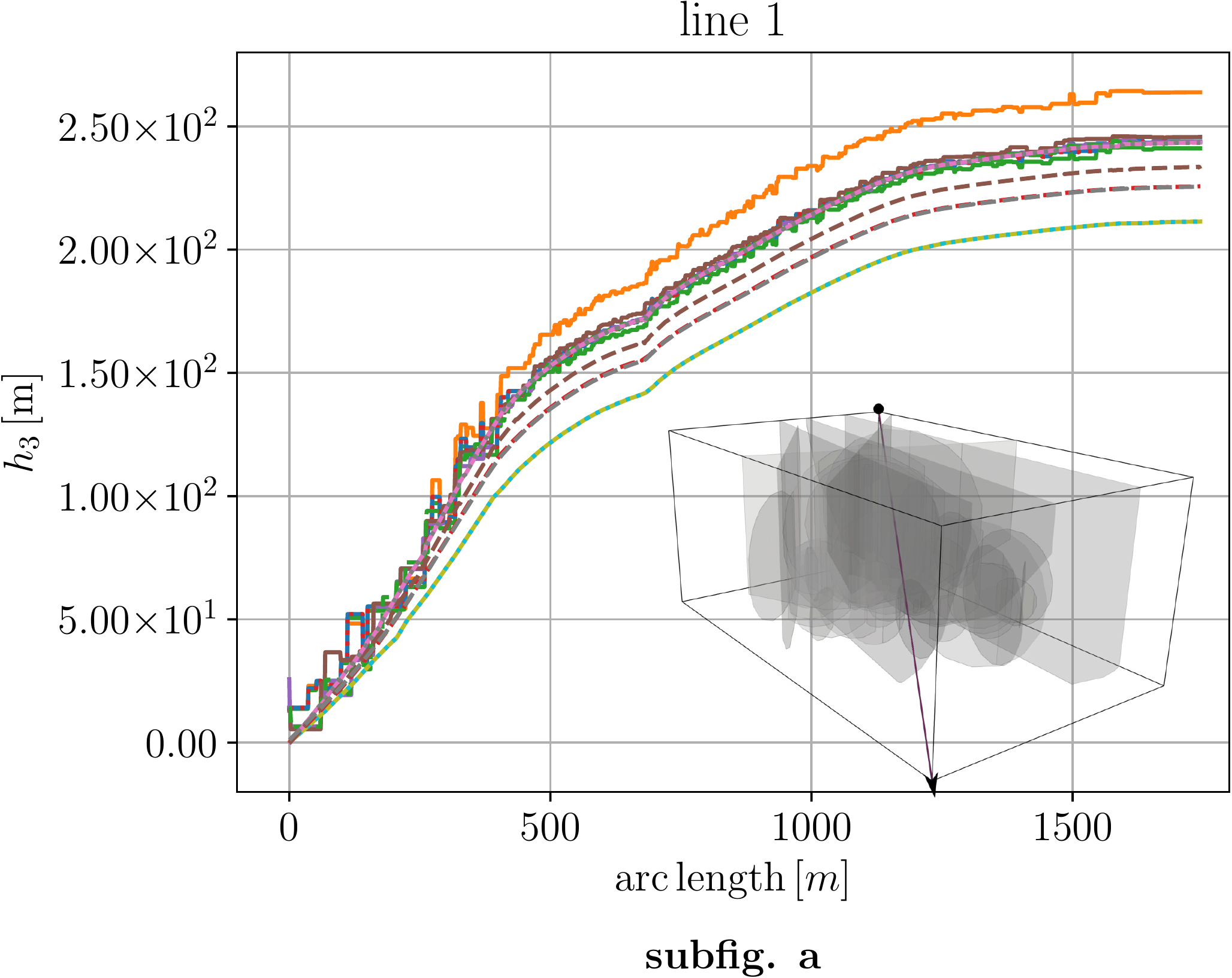}
\includegraphics[width=.45\textwidth]{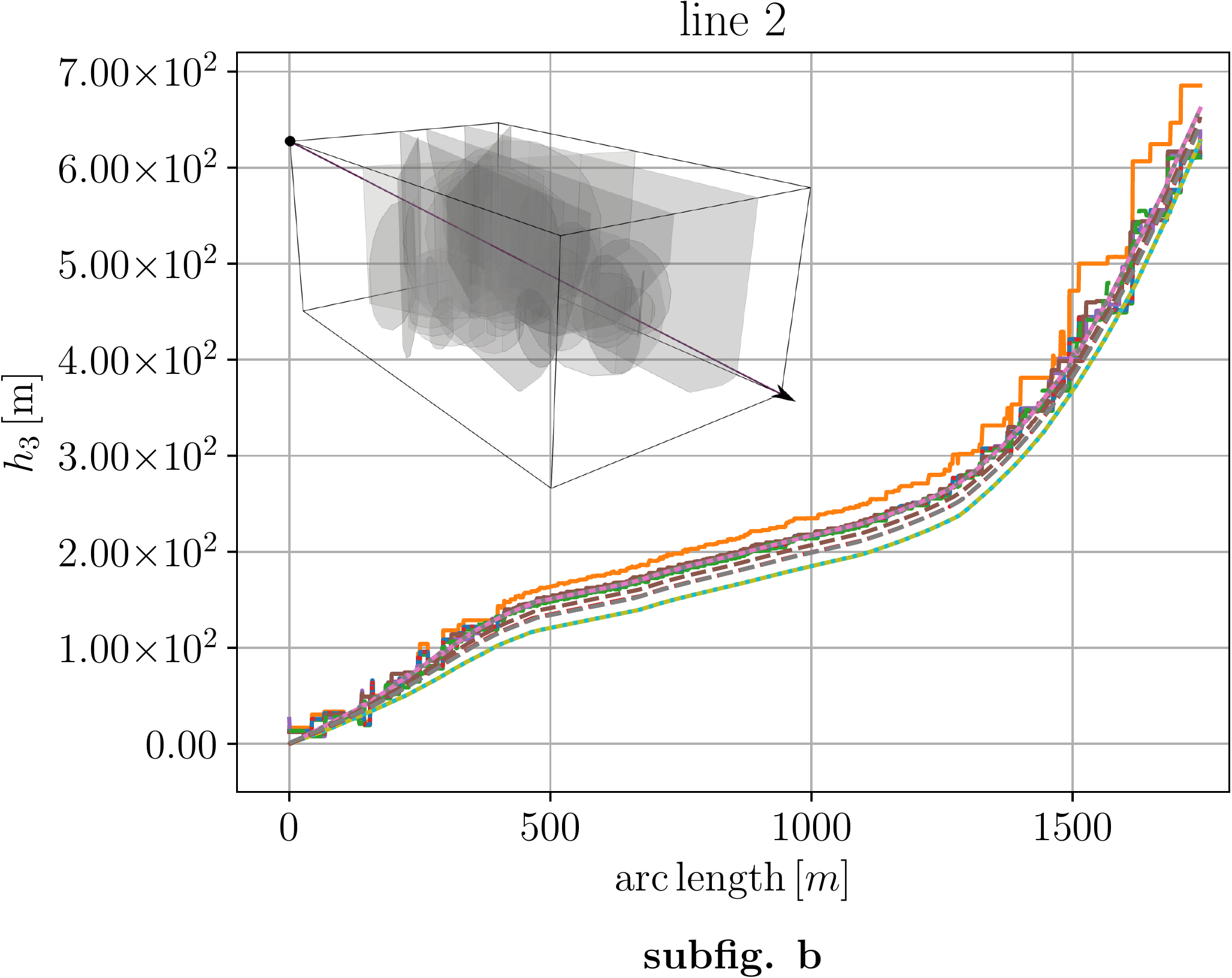}\\
\includegraphics[width=.45\textwidth]{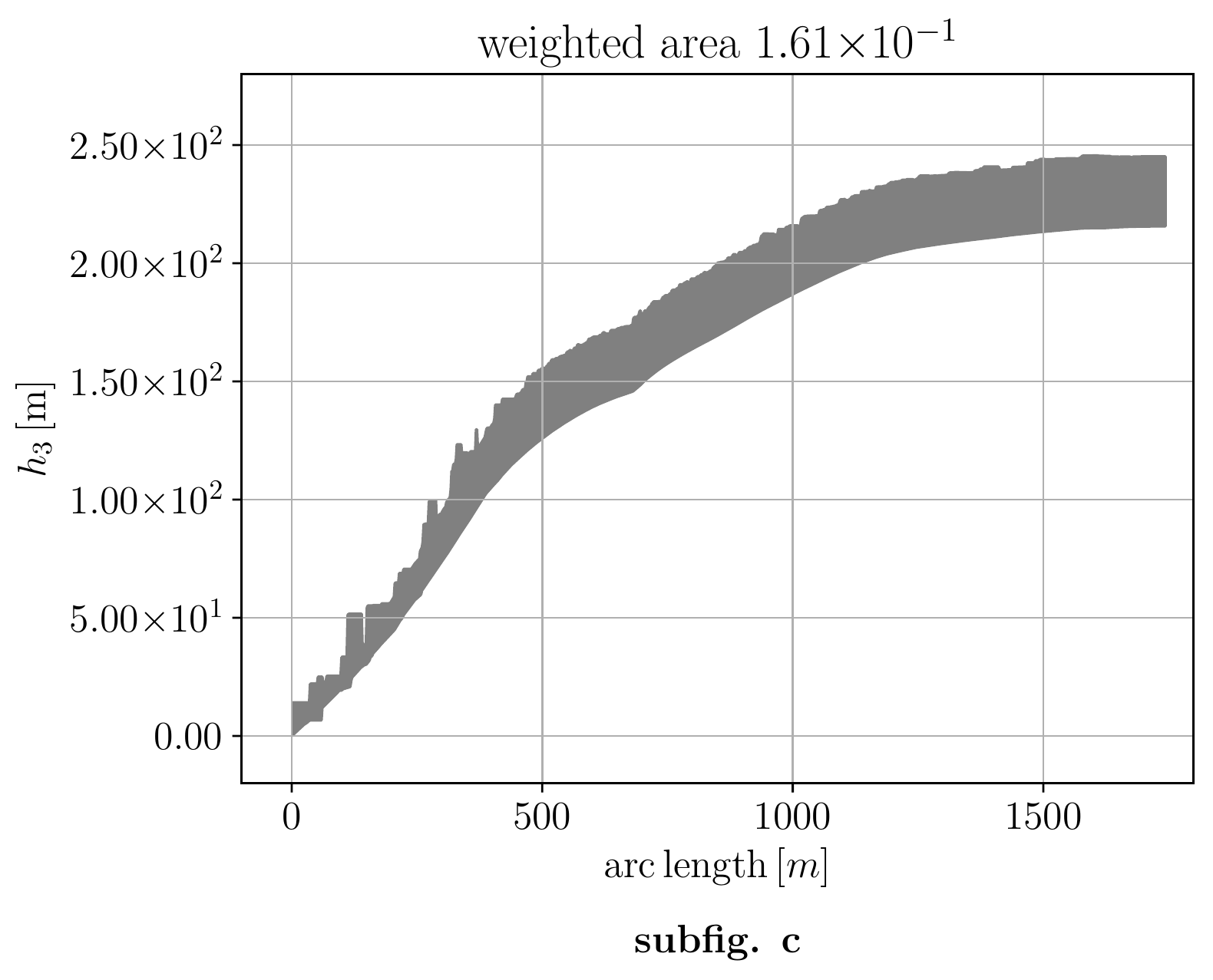}
\includegraphics[width=.45\textwidth]{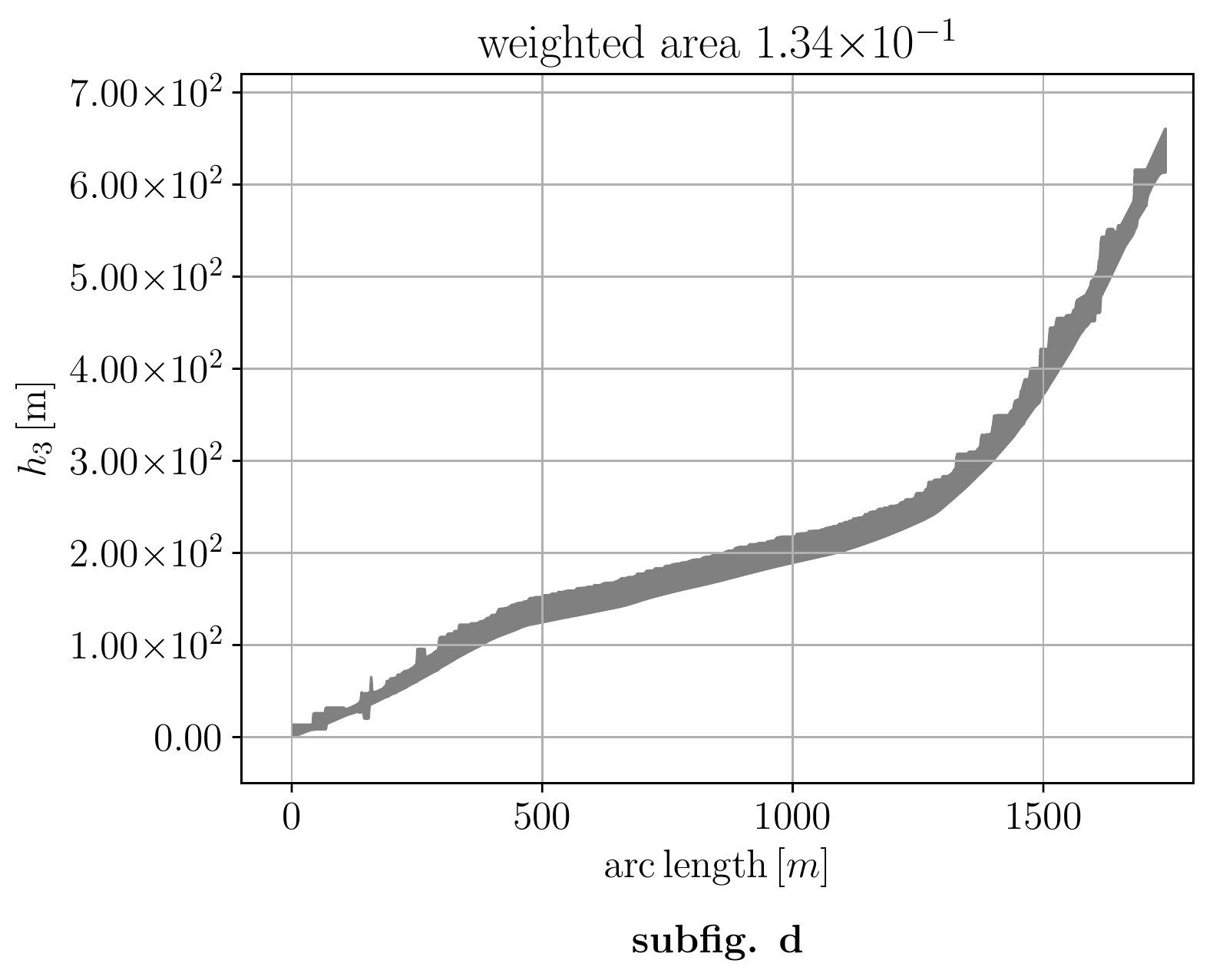}
\caption{Case 4 of Subsection \ref{subsec:field_network}. Hydraulic head profiles across the domain. Left: Profile from outlet $\partial\globalDomain_{out, 0}$ towards the opposite corner. Right: Profile from outlet $\partial\globalDomain_{out, 1}$ towards $\partial\globalDomain_{in}$.  Results of Subsection \ref{case4:a}.}
\label{fig:case4_pol}
\end{figure}

\paragraph{Hydraulic Head}\label{case4:a}
Figure \ref{fig:case4_pol} shows the hydraulic head along the two specified lines, together with the spread of the reported results.
Both lines start in points at the outflow boundaries where the hydraulic head is set to 0; the first line ends far away from the inlet, while the second ends at the inlet boundary. 
For the first line there are noticeable deviations for some of the solutions: 
The \UibTpfa scheme predicts a significantly higher hydraulic head drop, likely caused by the inconsistency of the scheme.
Conversely, the \UniceVagD and \UniceVagC methods underestimate the drop in hydraulic head compared to the average of the reported results, while there is only minor disagreement among the other methods.
On the second line, the \UibTpfa scheme overestimates the drop in hydraulic head over the domain, while the other methods are in very good agreement.

\begin{figure}[tp]
\centering
\includegraphics[width=.95\textwidth]{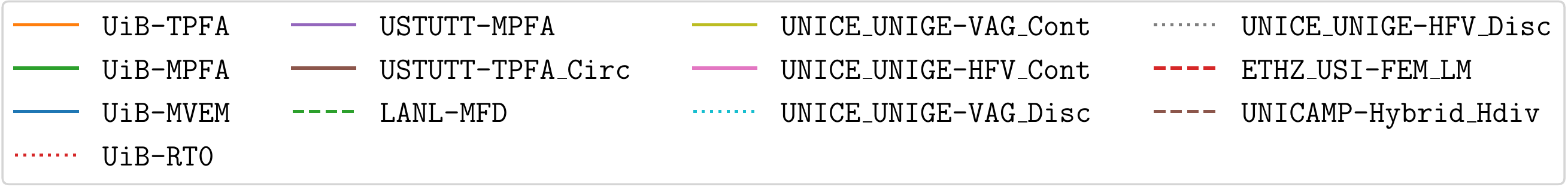}\\
\includegraphics[width=.95\textwidth]{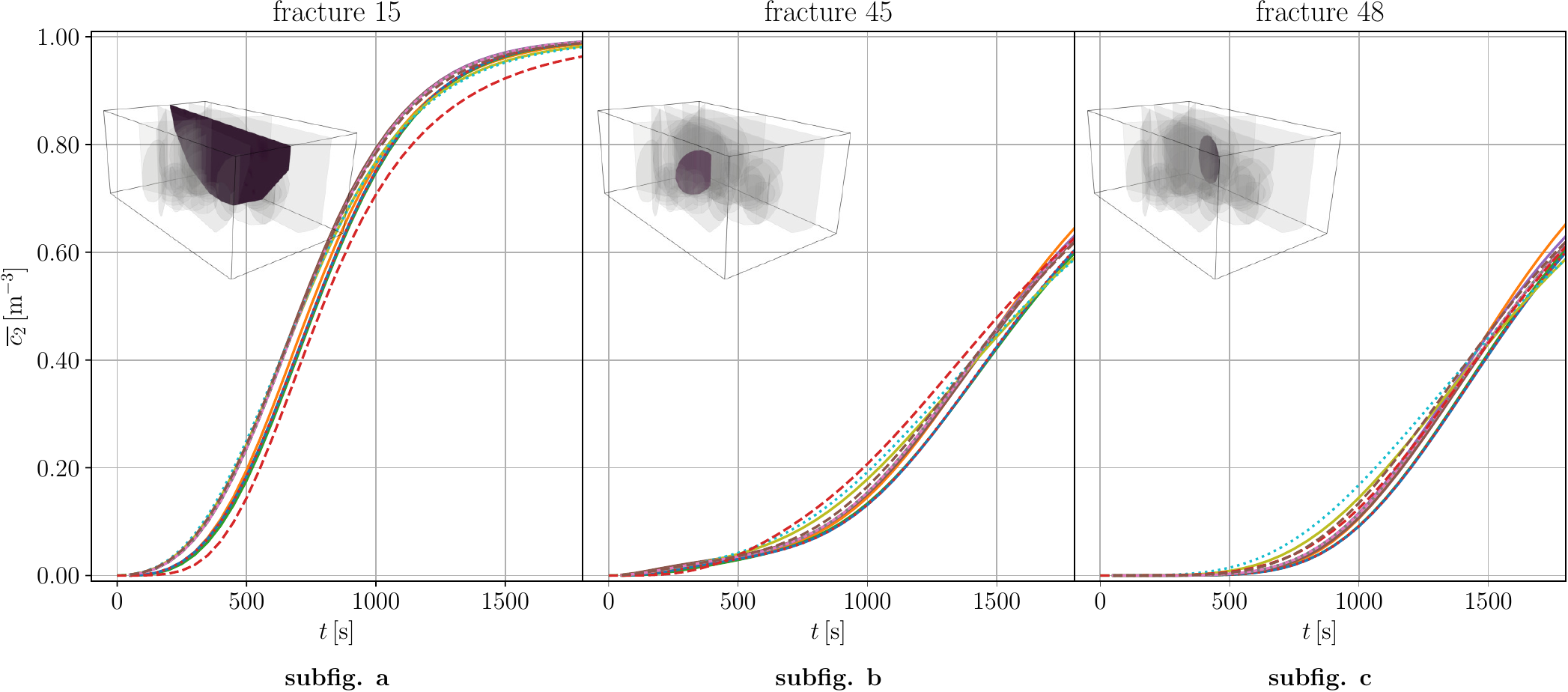}
\caption{Case 4 of Subsection \ref{subsec:field_network}. Mean concentration over time in three selected fractures with identification 15, 45, and 48.  Results of Subsection \ref{case4:b}.}
\label{fig:case4_pot}
\end{figure}

\paragraph{Concentration Plots}\label{case4:b}
The quality of the flux field is measured by the time series of average concentrations in the fracture planes, with good agreement among most of the methods.
Figure \ref{fig:case4_pot} shows the time evolution of concentration for three of the fractures, numbers 15, 45 and 48, which show the largest differences between the methods. 
The results produced by the \Ethz deviate slightly from the other methods on two of these figures, while \UniceVagD also shows a slight deviation for one of the figures.

\paragraph{Computational Cost}\label{case4:c}
Measures for the computational cost of the participating methods are given in Table \ref{tab:case4_comp_cost}.
Most of the groups used the provided mesh file. The \Unicamp method used a grid with only approximately 40\% of the cells in the provided grid. \Dtu employed almost seven times more 3d cells for its nodal-based method, yielding a number of degrees of freedom that is in the lower half with respect to all participating methods.
As in the previous test cases, there are significant differences in the number of unknowns and nonzero matrix elements among the methods.

\section{Summary of Results}\label{sec:discussion}
\begin{figure}[tp]
\centering
\includegraphics[width=.95\textwidth]{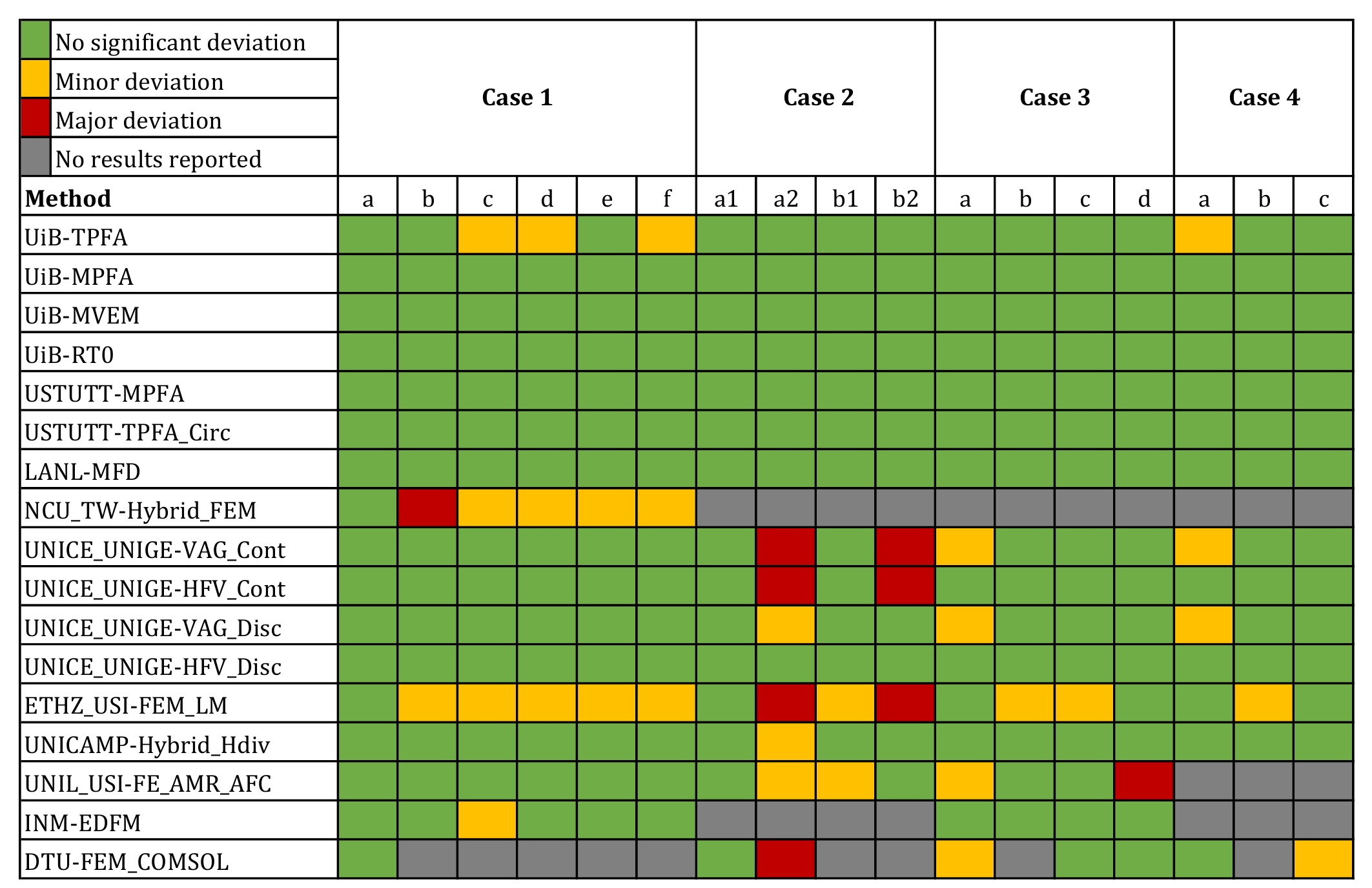}

\caption{Summary of the performance of all methods.}
\label{fig:summary_table}
\end{figure}

The performance of each method for all  test cases is indicated in Figure \ref{fig:summary_table}. We also list the main points emerging from the discussion of the results in Section \ref{sec:cases}:
\begin{enumerate}
\item Of the 17 schemes that participated in at least one of the test cases, 14 presented simulation results on all four cases.
\item Cases 3 and 4 pose the highest demands on the methods in terms of geometrical complexity.
Taken together, the cases point to the challenges inherent to DFM simulations and indicate the methods' robustness in this respect.
\item Not unexpectedly, fractures that act as barriers cause trouble for the methods that assume a continuous hydraulic head over the fracture, as seen in Case 2. Blocking fractures are outside the intended range of validity for these models, and alternative approaches should be sought for those cases.
\item Out of the 17 schemes, one 
is not mass conservative. There are no signs of the lack of conservation in the reported concentration fields, likely due to successful postprocessing of the flux fields. 
Nevertheless, for most of the test cases, the concentration fields reported by the nonconforming mesh method \Ethz deviate from the other reported results. 
\item The well-known inconsistency of the widely used two-point flux approximation is manifested in the underestimation of permeability in the hydraulic head results reported for \UibTpfa. The \StuttTpfa method circumvents this inconsistency by locating the hydraulic head values at the circumcenters of the tetrahedrons. However, this poses additional restrictions on the mesh.

\end{enumerate}

\section{Conclusion}\label{sec:conclusion}
This paper has presented a set of benchmark cases for the simulation of Darcy flow in three-dimensional fractured porous media.
The suite consists of one case with a single fracture, one case with 9 fractures and setups with conductive and blocking fractures, one case with 8 fractures designed to emphasize complex geometric details, and finally a case with 52 fractures, based on a real fracture network.
The metrics employed to measure discretization performance are (1) the profiles of the hydraulic head, (2) the quality of the flux field measured by simulation of passive tracers and (3) the computational cost as indicated by the number of degrees of freedom and matrix sparsity pattern. 
A total of 17 methods participated in the benchmark, spanning a wide range of discretization approaches for fractured media.
Although it was not possible to identify one approach as superior, the benchmark uncovered important differences between the methods.
The high number of participating methods and research groups proves that simulations in 3d media are fully feasible for a wide range of schemes and research codes.
For further development of discretization methods, 3d cases should therefore become a natural complement to the more traditional 2d simulation results. 

All data used in the benchmark can be found in the online repository \url{https://git.iws.uni-stuttgart.de/benchmarks/fracture-flow-3d.git}.
This includes the specification of benchmark case parameters, geometries, ready-made mesh generation files, applied metrics, and the results of all participating methods.
Additionally, Python scripts and Jupyter notebooks are provided which enable the reproduction of the result figures as well as the comparison of new computational results.
Therefore, we expect that the present work can serve as a reference point for the further development of discretization methods for fractured porous media.

\section{Acknowledgments}

I. Berre, A. Fumagalli, E. Keilegavlen, and I. Stefansson were supported by Norwegian Research Council grants 250223, 244129/E20 and 267908/E20.
The contribution of S. Burbulla was funded by the Deutsche Forschungsgemeinschaft (DFG, German Research Foundation) – Project Number 327154368 – SFB 1313.
P. Devloo was supported by Fapesp grant No 2017/15736-3, ANP/Petrobras grant No 2014/00090-2 and CNPq grant No 310369/2006-1.
O. Duran was supported by ANP/Petrobras grant No 2014/00090-2.
M. Favino was supported by the Swiss National Science Foundation (SNSF) grant PZ00P2\_180112.
M. Nestola and P.Zulian were supported by the SCCER-SoE and the Pasc Project FASTER.
C.F. Ni and I.H. Lee were partially funded by Institute of Nuclear Energy Research under grant NL1050288 and by Ministry of Science and Technology under grant MOST 108-2116-M-008-004- and MOST 106-2116-M-008-014-.
K. Nikitin and R. Yanbarisov were supported by Moscow Center for Fundamental and Applied Mathematics (agreement with the Ministry of Education and Science of the Russian Federation No. 075-15-2019-1624).
P. Schädle thanks the Werner Siemens Foundation for their endowment of the Geothermal Energy and Geofluids group at the Institute of Geophysics, ETH Zurich.

\section{Appendix}
\subsection{Measures of Computational Cost}\label{subsec:computational_cost}
This section provides three indicators related to computational cost: the number of cells (0d-3d), the number of degrees of freedom and the number of nonzero matrix entries. There is one table for each test case with data of all the participating methods at all refinement levels. For the equidimensional \Unil method, the cells listed as ``0d-2d cells'' are also three-dimensional cells that correspond to the fractures (``2d''), intersections of fractures (``1d'') and intersections of such intersections (``0d'').

\begin{table}
\begin{center}
\begin{tabular}{|l|l|l|l|l|l|l|l|}\hline
Method & Refinement & 0d cells & 1d cells & 2d cells & 3d cells & dofs & nnz \\ \hline
	\multirow{3}{*}{\UibTpfa}
			 & 0 & 0 & 0 & 112 & 1022 & 1358 & 6008\\
			 & 1 & 0 & 0 & 756 & 9438 & 11706 & 53904\\
			 & 2 & 0 & 0 & 4576 & 98311 & 112039 & 533547\\ \hline
	\multirow{3}{*}{\UibMpfa}
			 & 0 & 0 & 0 & 112 & 1022 & 1358 & 62200\\
			 & 1 & 0 & 0 & 756 & 9438 & 11706 & 672454\\
			 & 2 & 0 & 0 & 4576 & 98311 & 112039 & 7481237\\ \hline
    \multirow{3}{*}{\UibMVEM}
			 & 0 & 0 & 0 & 112 & 1022 & 3905 & 24435\\
			 & 1 & 0 & 0 & 756 & 9438 & 33651 & 222927\\
			 & 2 & 0 & 0 & 4576 & 98311 & 326561 & 2259630\\ \hline			 
	\multirow{3}{*}{\UibRT}
			 & 0 & 0 & 0 & 112 & 1022 & 3905 & 24435\\
			 & 1 & 0 & 0 & 756 & 9438 & 33651 & 222927\\
			 & 2 & 0 & 0 & 4576 & 98311 & 326561 & 2259623\\ \hline
	\multirow{3}{*}{\StuttMpfa}
			 & 0 & 0 & 0 & 100 & 1000 & 1100 & 22626\\
			 & 1 & 0 & 0 & 400 & 9600 & 10000 & 227354\\
			 & 2 & 0 & 0 & 3600 & 108000 & 111600 & 2731104\\ \hline
	\multirow{3}{*}{\StuttTpfa}
			 & 0 & 0 & 0 & 193 & 3400 & 3593 & 17373\\
			 & 1 & 0 & 0 & 448 & 9085 & 9533 & 46505\\
			 & 2 & 0 & 0 & 2582 & 104578 & 107160 & 530224\\ \hline
	\multirow{3}{*}{\Lanl}
			 & 0 & 0 & 0 & 100 & 1000 & 4400 & 51720\\
			 & 1 & 0 & 0 & 400 & 8000 & 34840 & 390840\\
			 & 2 & 0 & 0 & 1600 & 64000 & 267280 & 3035280\\ \hline
	\multirow{3}{*}{\Ncu}
			 & 0 & 0 & 0 & 625 & 9572 & 1840 & 25539\\
			 & 1 & 0 & 0 & 2453 & 65934 & 11537 & 169937\\
			 & 2 & 0 & 0 & 22262 & 638332 & 104581 & 1603776\\ \hline
	\multirow{3}{*}{\UniceVagC}
			 & 0 & 0 & 0 & 81 & 1134 & 1511 & 34085\\
			 & 1 & 0 & 0 & 361 & 10108 & 11721 & 288933\\
			 & 2 & 0 & 0 & 1849 & 103544 & 111233 & 2877105\\ \hline
	\multirow{3}{*}{\UniceHFVC}
			 & 0 & 0 & 0 & 81 & 1134 & 3870 & 39060\\
			 & 1 & 0 & 0 & 361 & 10108 & 32319 & 340879\\
			 & 2 & 0 & 0 & 1849 & 103544 & 320221 & 3454921\\ \hline
	\multirow{3}{*}{\UniceVagD}
			 & 0 & 0 & 0 & 81 & 1134 & 1943 & 43519\\
			 & 1 & 0 & 0 & 361 & 10108 & 13483 & 328867\\
			 & 2 & 0 & 0 & 1849 & 103544 & 119771 & 3073987\\ \hline
	\multirow{3}{*}{\UniceHFVD}
			 & 0 & 0 & 0 & 81 & 1134 & 4077 & 40041\\
			 & 1 & 0 & 0 & 361 & 10108 & 33231 & 345135\\
			 & 2 & 0 & 0 & 1849 & 103544 & 324779 & 3475475\\ \hline
	\multirow{3}{*}{\Ethz}
			 & 0 & 0 & 0 & 120 & 1000 & 1617 & 38834\\
			 & 1 & 0 & 0 & 480 & 10115 & 12714 & 335023\\
			 & 2 & 0 & 0 & 1920 & 93150 & 103470 & 2775270\\ \hline
	\multirow{3}{*}{\Unicamp}
			 & 0 & 0 & 0 & 526 & 1054 & 5968 & 114924\\
			 & 1 & 0 & 0 & 2884 & 10589 & 62164 & 1249536\\
			 & 2 & 0 & 0 & 15052 & 100273 & 604019 & 12448629\\ \hline
	\multirow{3}{*}{\Unil}
			 & 0 & 0 & 0 & 720 & 540 & 1857 & 49417\\
			 & 1 & 0 & 0 & 10880 & 38180 & 56947 & 1545935\\
			 & 2 & 0 & 0 & 39520 & 108671 & 579837 & 16878449\\ \hline
	\multirow{3}{*}{\Inm}
			 & 0 & 0 & 0 & 140 & 1000 & 1140 & 7666\\
			 & 1 & 0 & 0 & 720 & 10000 & 10720 & 73364\\
			 & 2 & 0 & 0 & 3800 & 100000 & 103800 & 719292\\ \hline
	\multirow{3}{*}{\Dtu}
			 & 0 & 0 & 0 & 0 & 1006 & 259 & 3082\\
			 & 1 & 0 & 0 & 0 & 10091 & 1931 & 26771\\
			 & 2 & 0 & 0 & 0 & 100014 & 17850 & 258202\\ \hline
\end{tabular}
\end{center}
\caption{Computational cost indicators for Case 1}
\label{tab:case1_comp_cost}
\end{table}

\begin{table}
\begin{center}
\begin{tabular}{|l|l|l|l|l|l|l|l|}\hline
Method & Refinement & 0d cells & 1d cells & 2d cells & 3d cells & dofs & nnz \\ \hline
	\multirow{3}{*}{\UibTpfa}
			 & 0 & 27 & 90 & 252 & 512 & 1820 & 8253\\
			 & 1 & 27 & 180 & 1008 & 4096 & 8074 & 43513\\
			 & 2 & 27 & 360 & 4032 & 32768 & 46622 & 281717\\ \hline 
	\multirow{3}{*}{\UibMpfa}
			 & 0 & 27 & 90 & 252 & 512 & 1820 & 8609\\
			 & 1 & 27 & 180 & 1008 & 4096 & 8074 & 44984\\
			 & 2 & 27 & 360 & 4032 & 32768 & 46622 & 287565\\ \hline 
	\multirow{3}{*}{\UibMVEM}
			 & 0 & 27 & 90 & 252 & 512 & 4706 & 20795\\
			 & 1 & 27 & 180 & 1008 & 4096 & 24862 & 118620\\
			 & 2 & 27 & 360 & 4032 & 32768 & 161414 & 806000\\ \hline 
	\multirow{3}{*}{\UibRT}
			 & 0 & 27 & 72 & 226 & 612 & 3970 & 21687\\
			 & 1 & 27 & 159 & 1192 & 5339 & 24727 & 153263\\
			 & 2 & 27 & 270 & 4536 & 39157 & 148245 & 980955\\ \hline 
	\multirow{3}{*}{\StuttMpfa}
			 & 0 & 0 & 0 & 284 & 843 & 1127 & 42060\\
			 & 1 & 0 & 0 & 686 & 3076 & 3762 & 207260\\
			 & 2 & 0 & 0 & 4578 & 38877 & 43455 & 2918322\\ \hline 
	\multirow{3}{*}{\StuttTpfa}
			 & 0 & 0 & 0 & 312 & 978 & 1290 & 7488\\
			 & 1 & 0 & 0 & 1206 & 4286 & 5492 & 31402\\
			 & 2 & 0 & 0 & 4578 & 38877 & 43455 & 226201\\ \hline
	\multirow{3}{*}{\Lanl}
			 & 0 & 0 & 0 & 434 & 628 & 2758 & 23246\\
			 & 1 & 0 & 0 & 1736 & 5024 & 18610 & 150314\\
			 & 2 & 0 & 0 & 6944 & 40192 & 134812 & 1062572\\ \hline
	\multirow{3}{*}{\UniceVagC}
			 & 0 & 0 & 0 & 252 & 512 & 974 & 22324\\
			 & 1 & 0 & 0 & 1008 & 4096 & 5902 & 143470\\
			 & 2 & 0 & 0 & 4032 & 32768 & 39908 & 1014088\\ \hline 
	\multirow{3}{*}{\UniceHFVC}
			 & 0 & 0 & 0 & 252 & 512 & 2223 & 22599\\
			 & 1 & 0 & 0 & 1008 & 4096 & 15048 & 157980\\
			 & 2 & 0 & 0 & 4032 & 32768 & 109368 & 1172592\\ \hline 
	\multirow{3}{*}{\UniceVagD}
			 & 0 & 0 & 0 & 252 & 512 & 2102 & 46348\\
			 & 1 & 0 & 0 & 1008 & 4096 & 10223 & 238891\\
			 & 2 & 0 & 0 & 4032 & 32768 & 56607 & 1390939\\ \hline 
	\multirow{3}{*}{\UniceHFVD}
			 & 0 & 0 & 0 & 252 & 512 & 2730 & 24138\\
			 & 1 & 0 & 0 & 1008 & 4096 & 17076 & 164148\\
			 & 2 & 0 & 0 & 4032 & 32768 & 117480 & 1197288\\ \hline 
	\multirow{3}{*}{\Ethz}
			 & 0 & 0 & 0 & 1212 & 512 & 3159 & 67183\\
			 & 1 & 0 & 0 & 1212 & 4096 & 7343 & 182793\\
			 & 2 & 0 & 0 & 1212 & 32768 & 38367 & 1036960\\ \hline  
	\multirow{3}{*}{\Unicamp}
			 & 0 & 27 & 69 & 534 & 923 & 6018 & 123312\\
			 & 1 & 27 & 90 & 1896 & 3912 & 23988 & 479322\\
			 & 2 & 27 & 249 & 10744 & 38742 & 236868 & 4830288\\ \hline 
	\multirow{3}{*}{\Unil}
			 & 0 & 1331 & 2787 & 6513 & 1745 & 16283 & 410491\\
			 & 1 & 1331 & 5211 & 20673 & 8129 & 45257 & 1180333\\
			 & 2 & 1331 & 10059 & 72033 & 47553 & 161805 & 4274281\\ \hline 
	\multirow{3}{*}{\Dtu}
			 & 0 & 0 & 0 & 0 & 550 & 129 & 1561\\
			 & 1 & 0 & 0 & 0 & 3881 & 836 & 10900\\
			 & 2 & 0 & 0 & 0 & 32147 & 6060 & 84954\\ \hline 
\end{tabular}
\end{center}
\caption{Computational cost indicators for Case 2}
\label{tab:case2_comp_cost}
\end{table}

\begin{table}[hbt]
\begin{center}
\begin{tabular}{|l|l|l|l|l|l|l|l|}\hline
Method & Refinement & 0d cells & 1d cells & 2d cells & 3d cells & dofs & nnz \\ \hline
	\multirow{2}{*}{\UibTpfa}
			 & 0 & 0 & 50 & 4305 & 31644 & 44786 & 207295\\
			 & 1 & 0 & 86 & 13731 & 138446 & 180024 & 849349\\ \hline
	\multirow{2}{*}{\UibMpfa}
			 & 0 & 0 & 50 & 4305 & 31644 & 44786 & 2596061\\
			 & 1 & 0 & 86 & 13731 & 138446 & 180024 & 11196843\\ \hline
	\multirow{2}{*}{\UibMVEM}
			 & 0 & 0 & 50 & 4305 & 31644 & 120696 & 818151\\
			 & 1 & 0 & 86 & 13731 & 138446 & 496032 & 3438098\\ \hline
	\multirow{2}{*}{\UibRT}
			 & 0 & 0 & 50 & 4305 & 31644 & 120696 & 818151\\
			 & 1 & 0 & 86 & 13731 & 138446 & 496032 & 3438098\\ \hline
	\multirow{2}{*}{\StuttMpfa}
			 & 0 & 0 & 0 & 4321 & 31942 & 36263 & 2459195\\
			 & 1 & 0 & 0 & 12147 & 131488 & 143635 & 10157331\\ \hline	 
	\multirow{2}{*}{\StuttTpfa}
			 & 0 & 0 & 0 & 4321 & 31942 & 36263 & 191147\\
			 & 1 & 0 & 0 & 12147 & 131488 & 143635 & 745375\\ \hline
	\multirow{2}{*}{\Lanl}
			 & 0 & 0 & 0 & 5617 & 21056 & 75878 & 607730\\
			 & 1 & 0 & 0 & 22468 & 168448 & 555887 & 4367379\\ \hline
	\multirow{2}{*}{\UniceVagC}
			 & 0 & 0 & 0 & 4321 & 31870 & 10213 & 130781\\
			 & 1 & 0 & 0 & 7711 & 150083 & 35485 & 479105\\ \hline	 
	\multirow{2}{*}{\UniceHFVC}
			 & 0 & 0 & 0 & 4321 & 31870 & 71708 & 504872\\
			 & 1 & 0 & 0 & 7711 & 150083 & 319175 & 2206691\\ \hline
	\multirow{2}{*}{\UniceVagD}
			 & 0 & 0 & 0 & 4321 & 31870 & 23302 & 400876\\
			 & 1 & 0 & 0 & 7711 & 150083 & 59187 & 966849\\ \hline
	\multirow{2}{*}{\UniceHFVD}
			 & 0 & 0 & 0 & 4321 & 31870 & 80538 & 532114\\
			 & 1 & 0 & 0 & 7711 & 150083 & 335599 & 2259971\\ \hline
	\multirow{2}{*}{\Ethz}
			 & 0 & 0 & 0 & 750 & 29295 & 33270 & 899809\\
			 & 1 & 0 & 0 & 3000 & 150930 & 163430 & 4421700\\ \hline
	\multirow{2}{*}{\Unicamp}
			 & 0 & 0 & 38 & 5580 & 24351 & 153519 & 3180847\\
			 & 1 & 0 & 51 & 23607 & 162773 & 994243 & 20600135\\ \hline
	\multirow{2}{*}{\Unil}
			 & 0 & 0 & 3877 & 323779 & 68386 & 86594 & 1206048\\
			 & 1 & 0 & 3877 & 323779 & 547088 & 148993 & 2202947\\ \hline
	\multirow{2}{*}{\Inm}
			 & 0 & 0 & 0 & 4036 & 29952 & 33988 & 240398\\
			 & 1 & 0 & 0 & 10732 & 149760 & 160492 & 1133364\\ \hline
	\multirow{2}{*}{\Dtu}
			 & 0 & 0 & 0 & 0 & 30984 & 5641 & 80669\\
			 & 1 & 0 & 0 & 0 & 150524 & 30379 & 469447\\ \hline
	\multirow{1}{*}{\texttt{USTUTT-MPFA-refined}}
 		     & 5 & 0 & 0 & 49428 & 980212 & 1029640 & 75207825\\ \hline 
\end{tabular}
\end{center}
\caption{Computational cost indicators for Case 3}
\label{tab:case3_comp_cost}
\end{table}

\begin{table}[hbt]
\begin{center}
\begin{tabular}{|l|l|l|l|l|l|l|l|}\hline
Method & 0d cells & 1d cells & 2d cells & 3d cells & dofs & nnz \\ \hline

	\UibTpfa   & 0 & 1601 & 52618 & 259409 & 424703 & 1950313\\
	\UibMpfa   & 0 & 1601 & 52618 & 259409 & 424703 & 22953336\\
	\UibMVEM   & 0 & 1601 & 52618 & 259409 & 1082740 & 7342691\\
	\UibRT     & 0 & 1601 & 52618 & 259409 & 1082740 & 7342691\\
	\StuttMpfa & 0 & 0 & 52618 & 259420 & 312038 & 21227071\\
	\StuttTpfa & 0 & 0 & 52618 & 259420 & 312038 & 1721932\\
	\Lanl      & 0 & 0 & 52070 & 260417 & 783158 & 7953396\\
	\UniceVagC & 0 & 0 & 52070 & 260431 & 95930 & 1237714\\
	\UniceHFVC & 0 & 0 & 52070 & 260431 & 600561 & 4349901\\
	\UniceVagD & 0 & 0 & 52070 & 260431 & 252326 & 4497980\\
	\UniceHFVD & 0 & 0 & 52070 & 260431 & 704813 & 4663105\\
	\Ethz      & 0 & 0 & 52618 & 212040 & 223532 & 5817930\\
	\Unicamp   & 0 & 938 & 24853 & 94294 & 629065 & 13233581\\
	\Dtu       & 0 & 0& 0 & 1860063 & 319489 & 4709565 \\

\hline
\end{tabular}
\end{center}
\caption{Computational cost indicators for Case 4}
\label{tab:case4_comp_cost}
\end{table}

\bibliographystyle{elsarticle-num}
\bibliography{literature}







\end{document}